\documentclass{amsart}

\makeindex

\usepackage[all]{xy}
\usepackage{amssymb}

\usepackage{amsmath}

\usepackage{amsthm}

\usepackage{amscd}

\usepackage{amsfonts}

\usepackage{amssymb}

\usepackage{mathrsfs}

\usepackage{enumerate}

\usepackage{color,cite,graphicx}

\newcommand{\calA}{\mathcal{A}}

\newcommand{\calH}{\mathcal{H}}

\newcommand{\be}{\mathcal{\begin{equation}}}
\newcommand{\ee}{\mathcal{\end{equation}}}

\newcommand{\bbT}{\mathbb{T}}

\newcommand{\bbF}{\mathbb{F}}
\newcommand{\bbC}{\mathbb{C}}

\newcommand{\bbP}{\mathbb{P}}
\newcommand{\bbQ}{\mathbb{Q}}

\newcommand{\bbZ}{\mathbb{Z}}

\newcommand{\frakh}{\mathfrak{h}}

\newcommand{\SL}{\textup{SL}}

\newcommand{\GL}{\textup{GL}}

\newcommand{\PGL}{\textup{PGL}}

\newcommand{\Sp}{\textup{Sp}}

\newcommand{\Kum}{\textup{Kum}}

\newcommand{\rank}{\textup{rank}}

\newcommand{\la}{\langle}
\newcommand{\ra}{\rangle}

\newcommand{\Aut}{\textup{Aut}}

\newcommand{\He}{\textup{He}}
\newcommand{\he}{\mathfrak{h}}

\newcommand{\Aff}{\textup{Aff}}

\newcommand{\PSp}{\textup{PSp}}

\newcommand{\Ca}{\textup{Ca}}

\newtheorem{Lem}{Lemma}[section]

    \newtheorem{Prop}[Lem]{Proposition}

    \newtheorem{Thm}[Lem]{Theorem}

    \newtheorem{lemma}{Lemma}

\theoremstyle{definition}

    \newtheorem{Rem}[Lem]{Remark}

\newcommand{\Ps}{\mathbb{P}^}
\newcommand{\Z}{\mathbb{Z}}

\newcommand{\lra}{\longrightarrow}

\newcommand{\C}{\mathbb C}
\newcommand{\e}{\epsilon}

\begin{document}
\title{The Hesse pencil of plane cubic curves}
\author{Michela Artebani}
\address{Departamento de Matem\'atica, Universidad de Concepci\'on, Casilla 160-C, Concepci\'on, Chile }
\email{martebani@udec.cl} 
\thanks{The first author was supported in part by: PRIN 2005: \emph{Spazi di moduli e teoria di Lie}, Indam (GNSAGA) and NSERC Discovery Grant of Noriko Yui
at Queen's University, Canada.}
\author{Igor Dolgachev}
\address{Department of Mathematics, University of Michigan, 525 E. University Av., Ann Arbor, Mi, 49109}
\email{idolga@umich.edu}
\thanks{The second author was supported in part by NSF grant 0245203.}

%\date{\today}
\begin{abstract}
This is a survey of the classical geometry of the Hesse configuration of  12 lines in the projective plane related to inflection points of a plane cubic curve. We also study two K3 surfaces with Picard number 20 which arise naturally in connection with the configuration. 
\end{abstract}
\maketitle 
\pagestyle{myheadings} 
\section{Introduction} In this paper we discuss some old and new results about the widely known \emph{Hesse configuration} of 9 points and 12 lines in the projective plane $\bbP^2(k)$: each point lies on 4 lines and each line contains 3 points, giving an abstract configuration $(12_3,9_4)$. Through most of the paper we will assume that $k$ is the field of complex numbers $\bbC$ although the configuration can be defined over any field containing three cubic roots of unity. The Hesse configuration can be realized by the 9 inflection points of a nonsingular projective plane curve of degree 3. This discovery 
%fact that these points give rise to a configuration $(12_3,9_4)$ 
is attributed to Colin Maclaurin (1698-1746) (see \cite{Pascal}, p. 384), however the configuration is named after Otto Hesse who was the first  to study its properties in \cite{H1},\cite{H2}\footnote{Not to be confused with another Hesse's configuration $(12_4,16_3)$,  also related to plane cubic curves, see  \cite{Do}.}. In particular, he proved  that the nine inflection points of a plane cubic curve form one orbit with respect to the projective group of the plane and can be taken as common inflection points of a pencil of cubic curves generated by the curve and its Hessian curve. In appropriate projective coordinates the \emph{Hesse pencil} is given by the equation
\begin{equation}\label{hessepencil}
\lambda (x^3+y^3+z^3)+\mu xyz = 0.
\end{equation}
The pencil was classically known as the \emph{syzygetic pencil}\footnote{The term ``syzygy'' originates in astronomy where three planets are in syzygy when they are in line. Sylvester used this word to express a linear relation between a form and its covariant. We will see later that the pencil contains the Hesse covariant of each of its members.}  of cubic curves (see \cite{CL}, p. 230 or \cite{Enriques}, p. 274), the name attributed to L. Cremona. We do not know who is  responsible for renaming the pencil, but apparently the new terminology  is widely accepted in  modern literature (see, for example, \cite{BH}).  

Recently Hesse pencils have become popular among number-theorists in connection with computational problems in arithmetic of elliptic curves (see, for example, \cite{Smart}), and also among theoretical physicists, for example  in connection with  homological mirror symmetry for elliptic curves (see \cite{Z}). 
%The Hesse equation of the  elliptic curve is deduced from the equation of a certain ring associated to automorphisms in the Fukaya category of the mirror symplectic torus. Moreover its higher degree analogs$$x_0^n+x_1^n+\ldots +x_n^n+\lambda x_0x_1\cdots x_n = 0$$
%define Calabi-Yau $(n-2)-$manifolds. The pencil is used to construct the mirror Calabi-Yau manifold as the quotient of the member $X_t$ of the pencil by a subgroup of $(\bbZ/n\bbZ)^{n-2}$ acting symplectically on $X_t$.

\begin{figure}\label{fig}
\begin{center}
\includegraphics[scale=.7]{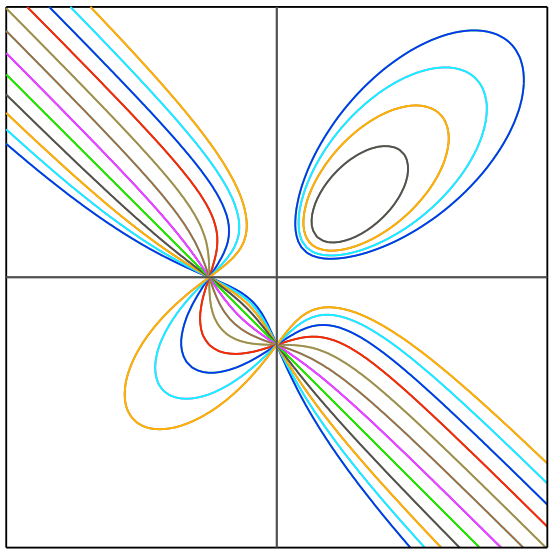}
%\epsffile[72 72 453 216]{he.eps}
\end{center}
\caption{The Hesse pencil} 
\end{figure}

The group of projective automorphisms which transform the Hesse pencil to itself is a group  $G_{216}$ of order 216 isomorphic to the group of affine transformations with determinant 1 of the projective plane over the field $\bbF_3$. 
This group was discovered by C. Jordan in 1878 \cite{J} who called it the \emph{Hessian group}. Its invariants were described in 1889 by H. Maschke in \cite{Maschke}. A detailed historical account and the first figure of the Hesse pencil can be found in \cite{G}. 

The projective action of the Hessian group comes from a linear action of a complex reflection group $\overline{G}_{216}$ of order 648 (n.25 in the Shephard-Todd list \cite{Shephard}) whose set of reflection hyperplanes consists of  the $12$ inflection lines of the Hesse configuration. The algebra of invariant polynomials of the group $\overline{G}_{216}$ is freely generated by three polynomials of degrees 6, 9, 12 (see \cite{Maschke}). An invariant polynomial of degree 6 defines a nonsingular plane sextic curve $C_6$. The double cover of the plane branched along the sextic curve $C_6$ is a K3 surface $X$ on which $G_{216}$ acts as a group of  automorphisms. Its subgroup $H = \bbF_3^2\rtimes Q_8$, where $Q_8$ is the Sylow 2-subgroup of $\SL(2,\bbF_3)$ isomorphic to the quaternion group of order 8, acts on the surface as a group of symplectic automorphisms. In fact, the group $\bbF_3^2\rtimes Q_8$ can be found in Mukai's list  \cite{M} of  finite groups which can be realized as   maximal finite groups of  symplectic automorphisms  of a complex K3 surface. 

The linear system of plane sextics with double points at 8 inflection points of a plane cubic is of projective dimension 3. The stabilizer $H$ of the ninth remaining inflection point in $G_{216}$ is isomorphic to $\SL(2,\bbF_3)$ and acts on this space by projective transformations. There is a unique invariant sextic $C'_6$ for this action, having cuspidal singularities at the inflection points. The double cover of the plane branched along $C'_6$ is birational to another  K3 surface $X'$ and the action of $H$ can be lifted to $X'$.  We show that $X'$ is birationally isomorphic to the quotient of $X$ by the subgroup $\bbF_3^2$ and the induced action of the quotient group $G_{216}/\bbF_3^2\cong \SL(2,\bbF_3)$ coincides with the action of $H$ on $X'$. Both K3 surfaces $X$ and $X'$  are singular in the sense of Shioda, i.e. the subgroup of algebraic cycles in the second cohomology group is of maximal possible rank equal to 20. We compute the intersection form defined by the cup-product on these subgroups. 

The invariant sextic $C_6$ cuts out a set of 18 points on each nonsingular member of the pencil. We explain its geometric meaning, a result which we were unable to find in the classical literature.\\
 
It is a pleasure to thank Bert van Geemen who was kind to provide us with his informal notes on this topic and to make many useful comments on our manuscript. We thank Noam Elkies and Matthias Schuett for their help in the proof of Theorem 7.10. We are also indebted to Thierry Vust for his numerous suggestions for improving the exposition of the paper. 

\section{The Hesse pencil}
Let $k$ be an algebraically closed field of characteristic different from 3 and $E$ be a nonsingular cubic in the projective plane $\bbP^2(k)$ defined by a homogeneous equation $F(x,y,z) = 0$ of degree 3.
The \emph{Hessian curve} $\He(E)$ of $E$ is the plane cubic curve defined by the equation $\He(F) = 0$, where $\He(F)$ is the determinant of the matrix of the second partial derivatives of $F$. The nine  points in $E\cap \He(E)$ are the \emph{inflection points} of $E$. 
 Fixing one of the inflection points $p_0$ defines a commutative group law $\oplus$ on $E$ with $p_0$ equal to the zero:
 $p\oplus q$ is the unique point $r$ such that $p_0, r$ and the third point of intersection in $\overline{p,q}\cap E$ lie on a line.
 It follows from this definition of the group law that each inflection point is a 3-torsion point and that the group $E[3]$  of 3-torsion points on $E$ is isomorphic to 
$(\bbZ/3\bbZ)^2$.  Any line $\overline{p,q}$ through two inflection points  intersects $E$ at another inflection point $r$ such that  $p,q,r$ form a coset  with respect to some subgroup of $E[3]$. Since we have 4 subgroups of order 3 in $(\bbZ/3\bbZ)^2$ we find 12 lines, each containing 3 inflection points. They are called the \emph{inflection lines} (or the \emph{Maclaurin lines} \cite{Enriques}) of $E$. Since each element in $(\bbZ/3\bbZ)^2$ is contained in 4 cosets, we see that  each inflection point is contained in four inflection lines. This gives the famous \emph{Hesse configuration} $(12_3,9_4)$ of 12 lines and 9 points in the projective plane.  It is easy to see that  this configuration is independent of the choice of the point $p_0$.\\

The \emph{Hesse pencil} is the one-dimensional linear system of plane cubic curves given by 
\begin{equation}\label{hessep}
E_{t_0,t_1}:\  t_0(x^3+y^3+z^3)+t_1xyz=0,\  \ (t_0,t_1)\in\Ps 1.\end{equation}
We use the affine parameter $\lambda = t_1/t_0$ and denote $E_{1,\lambda}$ by $E_\lambda$; the curve $xyz = 0$ is denoted by $E_\infty$.
Since the pencil is generated by the Fermat cubic $E_0$ and its Hessian, its nine base points are in the Hesse configuration. 
In fact, they are  the inflection points of any smooth curve in the pencil. In coordinates they are:
$$p_0=(0,1,-1),\ p_1=(0,1,-\e),\ p_2= (0,1,-\e^2), $$
$$p_3=(1,0,-1),\  p_4=(1,0,-\e^2),\  p_5=(1,0,-\e), $$
$$p_6=(1,-1,0),\  p_7=(1,-\e,0),\ p_8=(1,-\e^2,0),$$
where $\e$ denotes a primitive third root of 1.

If we fix the group law by choosing the point $p_0$ to be the zero point, then the set of inflection points is the group of 3-torsion points of each member of the Hesse pencil. Hence we can define an isomorphism
$$\xymatrix{\alpha:E_{\lambda}[3]_{p_0}\ar[r]& (\Z/3\Z)^2}$$
by sending the point $p_1 = (0,1,-\e)$ to $(1,0)$ and the point $p_3 = (1,0,-1)$ to $(0,1)$.
Under this isomorphism we can identify the nine base points with elements of $(\Z/3\Z)^2$ as follows
\begin{equation}\label{matrix2}
\begin{matrix}p_0&p_1&p_2\\
p_3&p_4&p_5\\
p_6&p_7&p_8\end{matrix} \  = \ 
\begin{matrix}(0,0)&(1,0)&(2,0)\\
(0,1)&(1,1)&(2,1)\\
(0,2)&(1,2)&(2,2).\end{matrix}
\end{equation}
It is now easy to see that any triple of base points which represents a row, a column, or a term in the expansion of the determinant of the matrix (\ref{matrix2}) spans an inflection line (cf. \cite{Miller}, p. 335).

The existence of an isomorphism $\alpha$ not depending on the member of the pencil can be interpreted by saying that the Hesse pencil is a family of elliptic curves together with a 3-level structure (i.e. a basis in the subgroup of 3-torsion points). In fact, in the following lemma we will prove that 
 any smooth plane cubic is projectively isomorphic to a member of the Hesse pencil.
It follows (see \cite{BH}) that its parameter space can  be naturally identified with  a smooth compactification of the fine moduli space $\calA_1(3)$ of elliptic curves with a 3-level  structure  (when $k = \bbC$  this is the \emph{modular curve} $X(3)$ \emph{of principal level} 3). 
 
\begin{lemma}
Any smooth plane cubic is projectively equivalent to a member of the Hesse pencil, i.e. it admits a Hesse canonical form\,\footnote{called the \emph{second canonical form} in \cite{Semple}, the first one being the Weierstrass form.}:
$$ x^3+y^3+z^3+\lambda xyz = 0. $$
\end{lemma}
\begin{proof}
We will follow the arguments from \cite{W}. Let $E$ be a nonsingular plane cubic. Given two inflection tangent lines for $E$ we can choose projective coordinates such that their equations are $x= 0$ and $y=0$. Then it is easy to see that the equation of $E$ can be written in the form
\begin{equation}\label{web}
F(x,y,z) = xy(ax+by+cz)+dz^3 = 0,
\end{equation}
 where $ax+by+cz = 0$ is a third inflection tangent line. Suppose  $c= 0$, then $ab\ne 0$ since otherwise the curve would be singular. Since a binary form of degree 3 with no multiple roots can be reduced, by a linear change of variables, to the form $x^3+y^3$,   the equation takes the form $x^3+y^3+dz^3 = 0$. After scaling the coordinate $z$, we arrive at a Hesse equation. So we may assume $c\ne 0$ and, after scaling the coordinate $z$,  that $c = 3$.
Let $\e$ be a primitive 3rd root of unity and define new coordinates $u,v$ by the formulae
$$ax+z = \e u+\e^2v, \quad by+z = \e^2u+\e v.$$
Then 
$$abF(x,y,z) = (\e u+\e^2v-z)(\e^2u+\e v-z)(-u-v +z)+dz^3$$
$$= -u^3-v^3 +(d+1)z^3-3uvz = 0.$$
Since the curve is nonsingular we have  $d\ne -1$. Therefore, after scaling the coordinate $z$, we get a Hesse equation for $E$:
\begin{equation}\label{hesse}
x^3+y^3+z^3+\lambda xyz = 0. \end{equation}\end{proof} 
%\noindent The equation (\ref{hesse}) is called \emph{Hesse canonical form}\footnote{called the \emph{second canonical form} in \cite{Semple}, the first one being the Weierstrass form.} of the plane cubic curve $E$.\\
%The curve given by this equation is singular if and only if $\lambda^3 = -27$.
%The \emph{Hesse pencil} is the one-dimensional linear system of cubic curves given by 
%\begin{equation}\label{hessep}
%E_{t_0,t_1}:\  t_0(x^3+y^3+z^3)+t_1xyz=0,\  \ (t_0,t_1)\in\Ps 1.\end{equation}
%We use the affine parameter $\lambda = t_1/t_0$ and denote $E_{1,\lambda}$ by $E_\lambda$. The curve $xyz = 0$ is denoted by $E_\infty$. 
%By the above, any nonsingular plane cubic is isomorphic to a member of this pencil. 
Assume additionally that the characteristic of the field $k$ is not equal to 2 or 3. Recall that a plane nonsingular cubic also admits the Weierstrass canonical form 
$$y^2z= x^3+axz^2+bz^3, \ \ 4a^3+27b^2 \ne 0.$$
%The absolute invariant
%$$j = \frac{4a^3}{4a^3+27b^2}$$
%distinguishes the projective equivalence classes of cubic curves. 
Projecting from the point $p_0$ we exhibit each curve of the Hesse pencil as a double cover of $\bbP^1$ branched at 4 points. By a standard procedure, this allows one to compute the Weierstrass form of any curve from the Hesse pencil:
\begin{equation}\label{weir}
y^2z= x^3+A(t_0,t_1)xz^2+B(t_0,t_1)z^3,
\end{equation}
where 
\begin{eqnarray}\label{binary}
A(t_0,t_1) &= &12u_1(u_0^3-u_1^3), \\
 B(t_0,t_1) &= & 2(u_0^6-20u_0^3u_1^3-8u_1^6), 
 \end{eqnarray}
and $ (t_0,t_1) =  (u_0,6u_1)$. The discriminant of the cubic curve given by \eqref{weir} is
$$\Delta = 4A^3+27B^2 = 2^23^3u_0^3(u_0^3+8u_1^3)^3,$$
its zeroes describe the singular members of the pencil.
The zeroes of the binary form $A(t_0,t_1)$ define the curves from the Hesse pencil which admit an automorphism of order 6 with a fixed point (\emph{equianharmonic cubics}).  For example, the Fermat curve 
$E_0: x^3+y^3+z^3 = 0$ is one of them. The zeroes of the binary form $B(t_0,t_1)$ define the curves from the Hesse pencil which admit an automorphism of order 4 with a fixed point (\emph{harmonic cubics}).   
The map 
$$j:\xymatrix{\bbP^1\ar[r]&\bbP^1}, \ (t_0,t_1) \mapsto (4A^3, 4A^3+27B^2)$$
coincides  (up to a scalar factor) with the map assigning to the elliptic curve $E_\lambda$ its $j$-invariant, which distinguishes the projective equivalence classes of cubic curves. \\

The Hesse pencil naturally defines a rational map
$$\xymatrix{\bbP^2 \ar@{-->}[r]&\bbP^1},\ (x,y,z) \mapsto (xyz, x^3+y^3+z^3)$$
which is not defined at the nine base points. Let 
$$\xymatrix{\pi:S(3) \ar[r]& \bbP^2}$$ 
be the blowing up of the base points. This is a  rational  surface such that the composition of rational maps $\xymatrix{S(3)\ar[r]&\bbP^2 \ar@{-->}[r]&\bbP^1}$ is a regular map
\begin{equation}\label{s3}
\xymatrix{\phi:S(3)\ar[r]&\bbP^1}\end{equation}
whose fibres are isomorphic to the members of the Hesse pencil. The map $\phi$ defines a structure of a minimal elliptic surface on $S(3)$.  Here and later we refer to  \cite{BPV}, \cite{Friedman},  \cite{Miranda}  or \cite{CD} for the theory of elliptic fibrations on algebraic surfaces. The surface 
$S(3)$ is a special case of an \emph{elliptic modular surface} $S(n)$ of level $n$ (see \cite{BH}, \cite{Shioda2}), isomorphic to  the universal family of elliptic curves with an  $n$-level.

There are four singular members in the Hesse pencil, each one is the  union of three lines:
$$\begin{array}{llc}
E_\infty:&\ \ xyz=0,\\
E_{-3}:&\ \ (x+y+z)(x+\e y+\e^2z)(x+\e^2y+\e z)=0,\\
E_{-3\e}:&\ \ (x+\e y+z)(x+\e^2 y+\e^2 z)(x+y+\e z)=0,\\
E_{-3\e^2}:&\ \ (x+\e^2 y+z)(x+\e y+\e z)(x+y+\e^2 z)=0.\\
\end{array}$$
We will call these singular members  the \emph{triangles} and denote them by  $T_1,\dots,T_4$, respectively.  The singular points of the triangles will be called  the \emph{vertices} of the triangles. They  are 
\begin{equation}\label{vertices}
\begin{array}{cccccc}
v_0=(1,0,0),& v_1=(0,1,0),& v_2=(0,0,1), \\
v_3=(1,1,1),& v_4=(1,\e,\e^2),& v_5=(1,\e^2,\e), \\
v_6=(\e,1,1),& v_7=(1,\e,1),& v_8=(1,1,\e), \\
v_9=(\e^2,1,1),& v_{10}=(1,\e^2,1),& v_{11}=(1,1,\e^2). 
\end{array}
\end{equation}
The 12 lines forming the triangles are the inflection lines of the Hesse configuration. If we fix a point $p_i$ as the origin in the group law of a nonsingular member of the pencil, then the side of a triangle $T_i$  passing through $p_i$ contains 3 base points forming a subgroup  of order 3, while the other sides of $T_i$ contain the cosets with respect to this subgroup.  
The triangles obviously give four singular fibres of Kodaira's type $I_3$ in the elliptic fibration $\phi$.

%Over complex numbers we can identify the base of the Hesse pencil with the modular curve $X(3) = \overline{H/\Gamma(3)}$ and the map $j$ coincides with the map $X(3)\to X(1) = \overline{H/\Gamma
%}$. The blow-up of the base points defines a rational elliptic surface $S(3)$ isomorphic to the modular elliptic surface of level 3 (see \cite{BH}).

\begin{Rem} The Hesse pencil makes sense over a field of any characteristic and is popular in number-theory and cryptography  for finding explicit algorithms  for computing the number of points of an elliptic curve over a finite field of characteristic 3 (see \cite{Ran}, \cite{Smart}). We are grateful to Kristian Ranestad for this comment. 

The proof of the existence of a Hesse equation for an elliptic curve $E$ over a  field of characteristic 3 goes through if we assume that $E$ is an ordinary elliptic curve with rational  3-torsion points. We find  equation \eqref{web} and check that it defines a nonsingular curve only if $abc\ne 0$. By scaling the variables we may assume that $a=b=-1, c = 1$. Next we use the variable change $z= 
u+x+y$ to transform the equation to the Hesse form 
$$xyu+d(u+x+y)^3 = xyu+d(u^3+x^3+y^3) = 0$$
The Hesse pencil \eqref{hessep} in characteristic 3 has two singular members: $(x+y+z)^3 = 0$ and $xyz = 0$. It has three base points $(1,-1,0), (0,1,-1),(1,0,-1)$, each of multiplicity 3,  which are the inflection points of all nonsingular members of the pencil. Blowing up the base points, including infinitely near base points, we get a rational elliptic surface. It has two singular fibres of Kodaira's types $IV^*$ and $I_2$. The fibre of type $IV^*$ has the invariant $\delta$ of wild ramification  equal to 1. This gives an example of a rational elliptic surface in characteristic 3 with finite Mordell-Weil group of sections (these surfaces are classified in \cite{Lang}). The  Mordell-Weil group  of our surface  is of order 3. 
\end{Rem}

%\begin{Rem} In modern time the Hesse pencil arises in homological mirror symmetry for elliptic curves (see \cite{Z}). The Hesse equation of the  elliptic curve is deduced from the equation of a certain ring associated to automorphisms in the Fukaya category of the mirror symplectic torus. Moreover its higher degree analogs
%$$x_0^n+x_1^n+\ldots +x_n^n+\lambda x_0x_1\cdots x_n = 0$$
%define Calabi-Yau $(n-2)-$manifolds. The pencil is used to construct the mirror Calabi-Yau manifold as the quotient of the member $X_t$ of the pencil by a subgroup of $(\bbZ/n\bbZ)^{n-2}$ acting symplectically on $X_t$.
%\end{Rem}

\section{The Hessian and the Cayleyan of a plane cubic} 
The first polar of a plane curve  $E$ with equation $F= 0$ with respect to a point $q = (a,b,c)\in \bbP^2$ is the curve
$P_q(E)$ defined by $aF_x'+bF_y'+cF_z' = 0$. It is easy to see that the Hessian curve $\He(E)$ of a plane cubic $E$ coincides with the locus of points $q$ such that the polar conic $P_q(E) = 0$ is reducible.
%is given by the equation $\He(F)= 0$, where $\He(F)$ is the determinant of the second partials of the polynomial $F$. 

If $E_\lambda$ is a member of the Hesse pencil, we find that $\He(E_\lambda)$ is the member $E_{\he(\lambda)}$ of the Hesse pencil, where
\begin{equation}\label{hes}
\he(\lambda) = -\frac{108+\lambda^3}{3\lambda^2}.
\end{equation}
 Let $p_i=(a,b,c)$ be one of the base points of the Hesse pencil. By computing the polar $P_{p_i}(E_\lambda)$ we find that it is equal to the union of the inflection tangent line $\bbT_{p_i}(E_\lambda)$ to the curve at the  point $p_i$ and the line $L_i: ax+by+cz = 0$. The  lines $L_0,\ldots, L_8$ are called the \emph{harmonic polars}.   
It follows easily from the known properties of the first polars (which can be checked directly in our case) that the line $L_i$ intersects the curve $E_{\lambda}$ at 3  points $q_j$ such that the tangent to the curve at $q_j$ contains $p_i$. Together with $p_i$ they form the group of 2-torsion points in the group law on the curve in which the origin is chosen to be the point $p_i$.

The harmonic polars, considered as points in the dual plane $\check{\bbP}^2$, give the set of base points of a  Hesse pencil in $\check{\bbP}^2$. Its inflection lines are the lines dual to the vertices of the inflection triangles given in \eqref{vertices}. If we identify the plane with its dual by means of the quadratic form $x^2+y^2+z^2$, the equation of the dual Hesse pencil coincides with the equation of the original pencil. For any nonsingular member of the Hesse pencil its nine tangents at the inflection points, considered as points in the dual plane, determine uniquely a member of the dual Hesse pencil. 

\begin{Rem} In the theory of line arrangements, the Hesse pencil defines two different arrangements (see \cite{BHH} and \cite{Hir}). The \emph{Hesse arrangement} consists of 12 lines (the inflection lines), it has 9 points of multiplicity 4 (the base points) and no other multiple points. The second arrangement is the dual of the Hesse arrangement, denoted by  $A_3^0(3)$.  It consists of 9 lines (the harmonic polars) and has 12 multiple points of multiplicity 3. Together these two arrangements form an abstract configuration $(12_3,9_4)$ which is a special case of a modular configuration (see \cite{Do}). In \cite{Hir} Hirzebruch constructs certain  finite covers of the plane with abelian Galois groups ramified over the lines of the Hesse configuration or its dual configuration. One of them, for each configuration, is a surface of general type  with universal cover isomorphic to a complex ball.
\end{Rem}

\begin{Prop} Let $E_\lambda$ be a nonsingular member of the Hesse pencil. Let $L_i\cap E_\lambda = \{q_1,q_2,q_3\}$ and $E_{\lambda_j},\ j = 1,2,3,$ be the curve from the Hesse pencil whose tangent at $p_i$ contains $q_j$. Then $\He(E_{\mu}) = E_\lambda$ if and only if $\mu\in\{\lambda_1,\lambda_2,\lambda_3\}$. 
\end{Prop}

\begin{proof} It is a straightforward computation. Because of the symmetry of the Hesse configuration, it is enough to consider the case when $i= 0$, i.e. $p_i = (0,1,-1)$.  We have that $L_0:y-z = 0$ and 
$L_0\cap E_\lambda$ is equal to the set of points $q_j = (1,y_j,y_j)$ satisfying
$1+2y_j^3+\lambda y_j^2 = 0$. The line $\overline{ p_0,q_j}$ has the equation $-2y_jx+y+z = 0$.
The curve $E_\mu$ from the Hesse pencil is tangent to this line at the point $(0,1,-1)$ if and only if 
$(-\mu,3,3) = (-2y_j,1,1)$, i.e. $y_j = \mu/6$. Thus 
$$\lambda = -\frac{1+2y_j^3}{y_j^2} = -\frac{108+\mu^3}{3\mu^2}.$$
Comparing with the formula \eqref{hes}, we see that $\he(\mu) = \lambda$. This proves the assertion.
\end{proof}

Let $E$ be a smooth plane cubic curve  which is not equianharmonic. Then $\He(E)$ is smooth and, for any $q\in \He(E)$, the polar conic $P_q(E)$ has one isolated singular point $s_q$.
In fact, $s_q$ lies on $\He(E)$ and the map $q\mapsto s_q$ is a fixed point free involution on $\He(E)$ (see, for example, \cite{DK}). If we fix a group law on $\He(E)$ with zero at $p_i$, then the map $q\mapsto s_q$ is the translation by a non-trivial 2-torsion point $\eta$. In the previous proposition this 2-torsion point is one of the intersection points of the harmonic polar  $L_i$ with $\He(E)$ such that $E$ is tangent to the line connecting this point with the inflection point $p_i$. 

The quotient $\He(E)/\la \eta\ra$ is isomorphic to the cubic curve in the dual plane $\check{\bbP}^2$ parametrizing the lines $\overline{ q,s_q}$. This curve is classically known as  the \emph{Cayleyan curve} of $E$. 
One can show that the Cayleyan curve also parametrizes the line components of reducible polar conics of $E$. In fact, the line $\overline{q,s_q}$ is a component of the polar conic $P_a(E)$, where $a$ is  the intersection point of the tangents of $\He(E)$ at $q$ and $s_q$. 
\begin{Prop}
If $E = E_\lambda$ is a member of the Hesse pencil, then its Cayleyan curve $\Ca(E_\lambda)$ is the member of the dual Hesse pencil corresponding to the parameter
\begin{equation}\label{cay}
\mathfrak{c}(\lambda) = \frac{54-\lambda^3}{9\lambda}.
\end{equation}
\end{Prop}
\proof
To see this,  following \cite{CL}, p. 245, we write the equation of the polar conic $P_q(E_{6\mu})$ with respect to a point $q=(u,v,w)$ 
$$u(x^2+2\mu yz)+v(y^2+2\mu xz)+w(z^2+2\mu xy) = 0.$$
It is a reducible conic if the equation decomposes into linear factors
$$u(x^2+2\mu yz)+v(y^2+2\mu xz)+w(z^2+2\mu xy) = (ax+by+cz)(\alpha x+\beta y+\gamma z).$$
This  happens if and only if 
$$\begin{pmatrix}u&2\mu w&2\mu v\\
2\mu w&v&2\mu u\\
2\mu v&2\mu u&w\end{pmatrix} = \begin{pmatrix}a\alpha&a\beta +b\alpha &a\gamma+c\alpha\\
a\beta+b\alpha&b\beta&c\beta+b\gamma\\
a\gamma+c\alpha&c\beta+b\gamma&c\gamma\end{pmatrix}.$$
Considering this as a system of linear equations in the variables $u,v,w,a,b,c$ we get the condition of solvability as the vanishing of the determinant
$$\begin{vmatrix}-1&0&0&\alpha&0&0\\
0&-1&0&0&\beta&0\\
0&0&-1&0&0&\gamma\\
-2\mu&0&0&0&\gamma&\beta\\
0&-2\mu&0&\gamma&0&\alpha\\
0&0&-2\mu&\beta&\alpha&0\end{vmatrix} = \mu(\alpha^3+\beta^3+\gamma^3)+(1-4\mu^3)\alpha\beta\gamma = 0.$$
If we take $(\alpha,\beta,\gamma)$ as the coordinates in the dual plane, this equation represents the equation of the Cayleyan curve because the line $\alpha x+\beta y+\gamma z$ is an irreducible component of a singular polar conic. Plugging in $\mu = \lambda/6$, we get \eqref{cay}.\qed\\

Note that the Cayleyan curve $\Ca(E_\lambda) = \He(E_\lambda)/\la \eta\ra$ comes with a distinguished nontrivial 2-torsion point, which is the image of the nontrivial coset of 2-torsion points on $\He(E_\lambda)$. This shows that $\Ca(E_\lambda) = \He(E'_\mu)$ for a uniquely defined member $E'_\mu$ of the dual Hesse pencil. The map $\alpha:\bbP^1\to \bbP^1, \lambda \mapsto \mu,$ gives an isomorphism between the spaces of parameters of the Hesse pencil and its dual pencil such that  
$\he(\alpha(\lambda)) = \mathfrak{c}(\lambda).$
One  checks that 
$$\he(-18/\lambda) = \mathfrak{c}(\lambda).$$

\begin{Rem} The Hesse pencil in the dual plane should not be confused with the (non-linear) pencil of the dual curves of members of the Hesse pencil. The dual curve of a nonsingular member $E_m = E_{m_0,3m_1}$ of the Hesse pencil is a plane curve of degree 6 with 9 cusps given by the equation
\begin{gather}
m_0^4(X_0^6+X_1^6+X_2^6)-m_0(2m_0^3+32m_1^3)(X_0^3X_1^3+X_0^3X_2^3+X_2^3X_1^3)\notag\\
-24m_0^2m_1^2X_0X_1X_2(X_0^3+X_1^3+X_2^3)-(24m_0^3m_1+48m_1^4)X_0^2X_1^2X_2^2 = 0.\label{dualcubic}
\end{gather}
This equation defines a surface $V$ in $\bbP^1\times \check{\bbP}^2$ of bi-degree $(4,6)$, the universal family of the pencil.  The projection to the first factor has fibres isomorphic to the dual curves of the members of the Hesse pencil, where the dual of a triangle becomes a triangle taken with multiplicity 2. The base points $p_i$  of the Hesse pencil define 9 lines $\ell_{p_i}$ in the dual plane and each of the 9 cusps of an irreducible member from \eqref{dualcubic} lies on one of these lines. The unique cubic passing through the nine cusps is the Cayleyan curve of the dual cubic. If $(m,x)\in V$, then the curve $E_m$ has the line $\ell_x$ (dual to $x$) as its tangent line. For a general point $x$, there will be 4 curves in the Hesse pencil tangent to this line, in fact the degree of the second projection $V\to \check{\bbP}^2$ is equal to 4. Each line $\ell_{p_i}$ lies in the branch locus of this map and its   preimage in $V$ has an irreducible  component $\bar{\ell}_{p_i}$ contained in the ramification locus. The surface $V$ is singular along the curves $\bar{\ell}_{p_i}$ and at the points corresponding to the vertices of the double triangles. One can show that a nonsingular minimal relative model of the elliptic surface $V\to \bbP^1$ is a rational elliptic surface isomorphic to $S(3)$. Thus, the dual of the Hesse pencil is the original Hesse pencil in disguise.
\end{Rem}

\begin{Rem} The iterations of the maps $\he:\bbP^1\to \bbP^1$ and $\mathfrak{c}:\bbP^1\to \bbP^1$ given by \eqref{hes} and \eqref{cay} were studied in \cite{Hollcroft}. They  give  interesting examples of  complex dynamics in one complex variable. The  critical points of $\he$ are the four equianharmonic cubics and its critical values correspond to the four triangles. Note that the set of triangles is invariant under this map. The set of critical points of $\mathfrak{c}$ is the set of triangles and it coincides with the set of critical values. The equianharmonic cubics are mapped to critical points. 
This shows that both maps are critically finite maps in the sense of Thurston (see \cite{McM}). 
\end{Rem}

\section{The Hessian group}\label{hessegroup}
The \emph{Hessian group} is the subgroup $G_{216}$ of   $\Aut(\Ps 2)\cong \PGL(3,\C)$ preserving the Hesse pencil.\footnote{Not to be confused with the Hesse group isomorphic to $\Sp(6,\bbF_2)$ and related to the 28 bitangents of a plane quartic.}  
The Hessian group acts on the space $\bbP^1$ of parameters of  the Hesse pencil, hence defines a homomorphism
\begin{equation}\label{al}
\xymatrix{\alpha:G_{216}\ar[r]& \Aut(\Ps 1)}.\end{equation}
Its kernel $K$ is generated by the transformations
\begin{eqnarray}\label{miss}
g_0(x,y,z)&=&(x, z, y),\\\notag
g_1(x,y,z)&=&(y,z,x), \\
g_2(x,y,z)&=&(x,\e y, \e^2 z).\notag
\end{eqnarray}
and  contains a normal subgroup of index 2
$$\Gamma=\langle g_1,g_2\rangle \cong (\Z/3\Z)^2.$$
If we use the group law with zero $p_0$ on a nonsingular member of the pencil, then $g_1$ induces the translation by the 3-torsion point $p_3$ and $g_2$ that by the point $p_1$.

%The   preimage $\overline{G}_{216}$ of $G_{216}$ in $\SL(3,\bbC)$ is a central extension $G_{216}$ of order $648$. It is a complex reflection group (see \cite{Springer}). The elements $g_1$ and $g_2$ generate a group $\Gamma$ of projective transformations isomorphic to $ (\Z/3\Z)^2$. Its   preimage in  $\overline{G}_{216}$  is isomorphic to the \emph{Heisenberg group} $\calH_3(3)$ of unipotent $3\times 3$ matrices over the field $\bbF_3$. 

The image of the homomorphism \eqref{al} is clearly contained in a finite subgroup of $\Aut(\bbP^1)$ isomorphic to the permutation group $S_4$. Note that it leaves the zeroes of the binary forms $A(t_0,t_1), B(t_0,t_1)$ from \eqref{binary} invariant. It is known that the group $S_4$ acts on $\bbP^1$ as an octahedral group, with orbits of cardinalities  $24, 12, 8, 6$, so it cannot leave invariant the zeroes of a binary form of degree 4. However, its subgroup $A_4$ acts as a tetrahedral group with orbits of cardinalities $12,6, 4,4$. This suggests that the image of \eqref{al} is indeed isomorphic to $A_4$.
In order to see it, it is enough to exhibit transformations from $G_{216}$ which are mapped to generators of $A_4$ of orders 2 and 3.
Here they are
 $$g_3=\left(\begin{array}{ccc}
 1&1&1\\
 1& \e & \e^2\\
 1& \e^2 & \e
 \end{array}\right),\ 
 g_4=\left(\begin{array}{ccc}
 1&0&0\\
0& \e & 0\\
 0& 0 & \e
 \end{array}\right).$$

 The group generated by $g_0,g_3,g_4$ is a central extension of degree two of $A_4$. It is isomorphic to the binary tetrahedral group and to the group $\SL(2, \bbF_3)$. Note that $g_3^2 = g_0$ so 
 $$G_{216} = \la g_1,g_2,g_3,g_4\ra.$$ It is clear that the order of $G_{216}$ is equal to the order of $K$ times the order of $A_4$  making it equal to $216$ that explains the notation.

\begin{Prop}
 The Hessian group $G_{216}$  is isomorphic to the semi-direct product
$$\Gamma \rtimes \SL(2, \bbF_3) $$
 %described by the  exact sequence
 %$$1\lra \Gamma\lra G_{216}\lra \SL(2, \bbF_3)\lra 1.$$
 where $\SL(2, \bbF_3)$ acts on $\Gamma\cong (\Z/3\Z)^2$ via the natural linear representation.
\end{Prop}

The Hessian group clearly acts on the set of nine points $p_i$,  giving a natural homomorphism from $G_{216}$ to $\Aff_2(3)$, the affine group of $\bbF_3^2$. In fact, the Hessian group is the subgroup of index 2 of $\Aff_2(3)$ of transformations with  linear part of determinant 1. In this action the group $G_{216}$ is realized as a 2-transitive subgroup of the permutation group $S_9$ on $\{0,1,\dots,8\}$ generated by permutations 
%$$T = (124)(568)(397)\ \mbox{ and }\ U = (456)(798)$$ 
$$T=(031)(475)(682)\ \mbox{ and }\  U=(147)(285)$$
%if we match the set of indices  $(0,1,2,\ldots,8)$ of the base points with the set $(1,4, 7, 2, 5, 8, 3,6,9)$ 
(see \cite{Coxeter}, 7.7). The stabilizer subgroup of the point $p_0$   is generated by $U$ and $TUT^{-1} = (354)(678)$ and coincides with  $\la g_3,g_4\ra$.
 \begin{Rem}  The group $\Aff_2(3)$ of order 432 that contains $G_{216}$ as a subgroup of order 2 is isomorphic to the Galois group of the equation of degree 9 defining the first  coordinates of the inflection points of a cubic with general coefficients in the affine plane \cite{Dickson}, \cite{W}. 
 \end{Rem}
%The transformations $g_3$ and $g_5=g_4^{-1}g_3 g_4$ generate a subgroup of $\PGL(2,\C)$ isomorphic to the quaternion group $Q_8$ with center generate by $g_0$. In fact, this is the $2$-Sylow subgroup of  $\SL(2, \bbF_3)$.

The Hessian group $G_{216}$, considered as a subgroup  $\PGL(3,\bbC)$, admits two different extensions to a subgroup of $\GL(3,\bbC)$ generated by complex reflections. The first group $\bar{G}_{216}$ is of order 648 and is generated by  reflections of order 3 (n.25 in Shephard-Todd's list \cite{Shephard}). The second group $\overline{G}_{216}'$ is of order 1296 and is generated by reflections of order 3 and  reflections of order 2 (n.26 in Shephard-Todd's list). The images of the reflection hyperplanes of $\overline{G}_{216}$ in the projective plane are the inflection lines, while the images of the reflection hyperplanes of $\overline{G}_{216}'$ are the inflection lines and the harmonic polars.

 The algebra of invariants of $\overline{G}_{216}$ is generated by three polynomials of degrees 6, 9 and 12 (see \cite{Maschke}, \cite{Springer})
$$\begin{array}{lll}
\Phi_6&= & x^6+y^6+z^6-10(x^3y^3+x^3z^3+y^3z^3),\\ \label{invariants}
\Phi_9&= &(x^3-y^3)(x^3-z^3)(y^3-z^3),\\ 
\Phi_{12}&= &(x^3+y^3+z^3)[(x^3+y^3+z^3)^3+216x^3y^3z^3].
\end{array}$$
Note that the curve $\Phi_9 = 0$ is the union of the nine harmonic polars  $L_i$ and the curve $\Phi_{12} = 0$ is the union of the four equianharmonic members of the pencil.  The union of the 12 inflection lines is obviously invariant with respect to $G_{216}$, however the corresponding polynomial $\Phi_{12}'$ of degree 12  is not an invariant but a relative invariant (i.e. the elements of $\overline{G}_{216}$ transform the polynomial to a constant multiple).

The algebra of invariants of the second complex reflection group $\overline{G}_{216}'$ is generated by $\Phi_6, \Phi_{12}'$ and a polynomial of degree 18
$$\Phi_{18} = (x^3+y^3+z^3)^6-540x^3y^3z^3(x^3+y^3+z^3)^3-5832x^6y^6z^6.$$
The curve $\Phi_{18} = 0$ is the union of the six harmonic cubics in the pencil. We will give later a geometric meaning to the 18 intersection points  of the curve defined by $\Phi_6 = 0$ with nonsingular members of the pencil.

A third natural linear extension of the group $G_{216}$ is the preimage $G_{216}'$ of the group under the projection $\SL(3,\bbC)\to \PGL(3,\bbC)$. This is a group of order 648 isomorphic to the central extension  $3G_{216}$ of $G_{216}$,  however it is not isomorphic to $\overline{G}_{216}$. The   preimage of the subgroup $\Gamma$ in $3G_{216}$ is a non-abelian group of order 27 isomorphic to the \emph{Heisenberg group} $\calH_3(3)$ of unipotent $3\times 3$-matrices with entries in $\bbF_3$. The group $G_{216}'$ is then isomorphic to the semi-direct product $\calH_3(3)\rtimes \SL(2,\bbF_3)$ and it is generated by  $g_1,g_2, \frac{1}{\e-\e^2}g_3,e^{2\pi i/9}g_4$ considered as  linear transformations.
\begin{Rem} Classical geometers were used to define a projective 
transformation as a pair  consisting of a nondegenerate quadric  in the projective space and  a  nondegenerate quadric in the dual projective space.  If $\bbP^n = \bbP(V)$, then the first quadric is given by a quadratic form on $V$ which defines a linear map $\phi:V\to V^*$. The second quadric defines a linear map $V^*\to V$ and the composition with the first one is a linear map $V\to V$.  In \cite{G} the Hessian group is given by a set of 36 conics  which are identified with conics in the dual plane $\check{\bbP}^2$ by means of an isomorphism $\bbP^2\to \check{\bbP}^2$ defined by the conic $x_0^2+x_1^2+x_2^2 = 0$. These conics are the polars of four equianharmonic cubics in the pencil with respect to the 12 vertices of the inflection  triangles. The 12 of them which are double lines have to be omitted.\end{Rem}

\smallskip
It is known that the simple group $G = \PSp(4,\bbF_3)$ of order 25,920 has two maximal subgroups of index 40. One of them is isomorphic to the complex reflection  group $\overline{G}_{216}$ of order 648. It has the following beautiful realization in terms of complex reflection groups in dimensions 4 and 5. 

It is known that the group $\bbZ/3\bbZ \times Ê\Sp(4,\bbF_3)$ is isomorphic to a complex reflection group in $\bbC^4$ Êwith 40 reflection hyperplanes of order 3 (n.32 in Shephard-Todd's list). This defines a projective representation of $G$ in $\bbP^3$ and the their stabilizer subgroups  of the reflection projective planes are isomorphic to $\overline{G}_{216}$. The reflection planes cut out on each fixed reflection plane the extended Hesse configuration of  12 inflection lines and 9 harmonic polars \cite{Maschke}, p. 334. 
 
Also is known that the group $\bbZ/2\bbZ \times G \not\cong \Sp(4,\bbF_3)$ is  isomorphic to a complex reflection group in $\bbC^5$ with 45 reflections hyperplanes of order 2 (n.33 in Shephard-Todd's list). This defines a projective representation of $G$ in $\bbP^4$. The algebra of invariant polynomials with respect to the complex reflection group $\bbZ/2\bbZ \times G$ was computed by  Burkhardt \cite{Burkhardt}. The smallest degree invariant is of degree 4. Its zero locus in $\bbP^4$ is the famous Burkhardt  quartic hypersurface with 45 nodes where 12 reflection hyperplanes meet. There are 40 planes forming one orbit, each containing 9 nodes. Each such plane contains 12 lines cut out by the reflection hyperplenes. They form the Hesse configuration with the 9 points equal to the set of base points of the Hesse pencil.

One can find an excellent exposition of the results of Maschke and Burkhardt in \cite{Hunt}. There is also a beautiful interpretation of the geometry of the two complex reflection groups in terms of the moduli space $\calA_2(3)$ of principally polarized abelian surfaces with some 3-level structure (see \cite{FSM},  \cite{Geer}). 
For example, one can identify $\calA_2(3)$ with  an open subset of the  Burkhardt quartic  whose complement is equal to the union of the 40 planes. 

\section{The quotient plane}
Consider the blowing up $\pi:S(3)\lra \Ps 2$
of  the base points $p_i$ of the Hesse pencil and the elliptic fibration \eqref{s3}
$$\xymatrix{\phi:S(3)\ar[r]& \Ps 1},\  \ (x,y,z)\mapsto (xyz, x^3+y^3+z^3).$$
 The action of the group $\Gamma$ on $\bbP^2$ lifts to an action on $S(3)$. Fixing one section of $\phi$ (i.e. one point $p_i$), the group $\Gamma$ is identified with the Mordell-Weil group of the elliptic surface and its action with the translation action.
Let 
$$\xymatrix{\bar\phi: S(3)/\Gamma\ar[r]& \Ps 1},\hspace{1cm} \xymatrix{\bar\pi:S(3)/\Gamma\ar[r]&\Ps 2/\Gamma},$$ 
be the morphisms induced by $\phi$ and $\pi$, respectively.

\begin{Prop}\label{quo2} The quotient surfaces  $\Ps 2/\Gamma$ and $S(3)/\Gamma$  have 4 singular points of type $A_2$ given by the orbits of the vertices (\ref{vertices}). The minimal resolution of singularities  are isomorphic to a Del Pezzo surface $S$ of degree 1 and to $S(3)$, respectively.
Up to these resolutions, $\bar\phi$  is isomorphic to $\phi$ and $\bar\pi$ is the blowing up of $S$ in one point, the $\Gamma$-orbit of the points $p_i$. 
\end{Prop}

\begin{proof} 
The group $\Gamma$ preserves each singular member of the Hesse pencil and any of its subgroups of order 3 leaves invariant the vertices of one of the triangles. Without loss of generality we may assume that the triangle is $xyz = 0$. Then the subgroup of $\Gamma$ stabilizing its vertices is generated by the transformation $g_2$, which acts locally at the point $y= z = 0$  by the formula $(y,z)\mapsto (\e y, \e^2z)$.
It follows that the orbits of the vertices give 4 singular points of type $A_2$ in $\Ps 2/\Gamma$ and $S(3)/\Gamma$, locally given by the equation $uv+w^3 = 0$.

 Let $E$ be an elliptic curve with a group law and $[n]:E\to E$ be the map $x\mapsto nx$. It is known that this map is a surjective map of algebraic groups with kernel equal to the group of $n$-torsion points. Its degree is equal to $n^2$  if $n$ is coprime to the characteristic.  In our case the quotient map by $\Gamma$ acts on each member of the Hesse pencil as the map $[3]$. This implies that the quotient of the surface $S(3)$ by the group $\Gamma$ is isomorphic to $S(3)$ over the open subset $U = \bbP^1\setminus \{\Delta(t_0,t_1) = 0\}$.

%Any singular fibre of $\phi:S(3)\to \bbP^1$ is invariant with respect to $\Gamma$ and there is subgroup of order 3 which leaves invariant its set of singular points.  Without loss of generality we may assume that the fibre originates from the triangle $xyz = 0$. Then the subgroup stabilizing the components of this  fibre is generated by the transformation $g_2$, which acts locally at the singular point $y= z = 0$  by the formula $(y,z)\mapsto (\e y, \e^2z)$. The orbit of this point in $S(3)/\Gamma$ has a singularity of type $A_2$, locally given by the equation $uv+w^3 = 0$. 
The map $\bar{\phi}:S(3)/\Gamma \to \bbP^1$ induced by the map $\phi$ has four  singular fibres. Each fibre is  an irreducible  rational curve with a double point which is  a singular point of the surface of type $A_2$.  Let $\sigma:S(3)'\to S(3)/\Gamma$ be a minimal resolution of the four singular points of $S(3)/\Gamma$.  The composition  $\bar{\phi}\circ \sigma:S(3)'\to \bbP^1$ is  an elliptic surface  isomorphic to $\phi:S(3)\to \bbP^1$ over the open subset $U$ of the base $\bbP^1$. Moreover, $\bar{\phi}\circ \sigma$ and $\phi$ have singular fibres of the same types, thus $S(3)'$ is a minimal elliptic surface. Since it is known that a birational isomorphism of minimal elliptic surfaces is an isomorphism, this implies that $\bar{\phi}\circ \sigma$ is isomorphic to $\phi$.

The minimal resolution $S$ of $\Ps 2 /\Gamma$ contains a pencil of cubic curves intersecting in one point $q_0$, the orbit of the points $p_i$. Hence it follows easily (see for example \cite{CD}) that $S$ is isomorphic to a Del Pezzo surface of degree one and $\bar\pi$ is the blowing up of the point $q_0$.\end{proof}

Let $\pi':S(3)'\to \bbP^2$ be the contraction of the 9 sections $E_0,\dots,E_8$ of the elliptic fibration $\bar{\phi}\circ \sigma$ to the points $q_0,\ldots,q_8$  in $\Ps 2$,  the base points of the Hesse pencil in the second copy of $\bbP^2$. 

By Proposition \ref{quo2} the following diagram is commutative:
\begin{equation}\label{diagram}
\xymatrix{S(3)\ar[d]_{\pi}\ar[r]^r&
S(3)/\Gamma \ar[d]^{\bar{\pi}}&S(3)'\ar[d]_\alpha\ar[l]_\sigma\ar[dr]^{\pi'}&\\
  \Ps 2 \ar[r]^p& \bbP^2/\Gamma &S\ar[l]_\beta\ar[r]^\gamma&\bbP^2.}\end{equation}
 %Here $\pi':S(3)'\to \bbP^2$ is the blow-up of the 9 base points $q_0,\ldots,q_8$ of the Hesse pencil in the second copy of $\bbP^2$, the image of the elliptic pencil on $S(3)'$. 
Here $p$ is the quotient map by $\Gamma$, $\beta$ is a minimal resolution of singularities of the orbit space $\bbP^2/\Gamma$, $\alpha $ is the blow-up of the point $q_0$ on $S$ and
 $\gamma$ is the blow-up of $q_1,\dots,q_8$ (see the notation in the proof of Proposition \ref{quo2}). %The map $\beta:S\to \bbP^2/\Gamma$ is a minimal resolution of singularities of the orbit space $\bbP^2/\Gamma$ and the map $\alpha:S(3)'\to S$ is the blow-up of one point on $S$ whose image under the map $\gamma$ is the ninth base point $q_0$ of the Hesse pencil, the $\Gamma$-orbit of the nine base points $p_0,\ldots,p_8$ of the original Hesse pencil. \\

%Fix a group law on nonsingular  members $F_t$ of the Hesse pencil with the origin $p_0$.  A 9-torsion point 
\begin{Prop}\label{hc} 
%Let  $ E_i = \pi'{}^{-1}(q_i),\  i = 0,\ldots,8,$ be the nine sections of the elliptic pencil on $S(3)'$. 
The curves $B_i = p^{-1}(\beta(\alpha(E_i)) ),\ i = 1,\ldots,8,$ are  plane cubic curves with equations
$$\begin{array}{ll}
B_1:\  x^3+\e y^3+\e^2 z^3 = 0,& B_5:\  x^3+\e^2 y^3+\e z^3 = 0,\\
B_2:\  x^2y+y^2z+z^2x = 0,& B_6: \  x^2z+y^2x+z^2y = 0,\\
B_3:\  x^2y+\e^2 y^2z+\e z^2x = 0, & B_7:\  x^2z+\e y^2x+\e^2 z^2y = 0,\\
B_4:\  x^2y+\e y^2z+\e^2 z^2x = 0,& B_8:\ x^2z+\e^2 y^2 x+\e z^2 y = 0.
\end{array}$$
The union of the eight cubics $B_i$ cuts out on each nonsingular member of the Hesse pencil the set of points of order 9 in the group law with the point $p_0$ as the origin.  

Each of them has one of the triangles of the Hesse pencil as inflection triangle and it is inscribed and circumscribed to the other three triangles (i.e. is tangent to one side of the triangle at each vertex). 
 \end{Prop}
 
\begin{proof} Recall that the sections $E_1,\ldots,E_8$ on $S(3)'$ are non-trivial 3-torsion sections (the zero section is equal to $E_0$). The   preimage $\bar B_i$  of $E_i$ under the map $r^{-1}\circ \sigma$ cuts out on each nonsingular fibre the $\Gamma$-orbit of a  point of order 9. Thus the image $B_i$ of $\bar B_i$ in $\bbP^2$ is a plane cubic cutting out the $\Gamma$-orbit of a point of order 9 on each nonsingular member of the Hesse pencil. 

Let $E$ be a  nonsingular member of  the Hesse pencil. Take a point $p\in E$ and let $q\ne p$ be the intersection of $E$ with the tangent line at $p$. Let $r\ne q$ be the intersection of $E$ with the tangent line at $q$. Finally, let $s\ne r$ be the intersection of $E$ with the tangent line at $r$.  
%Then $p$ is a point of order 9 if and only if $s=p$ (this explains why the classics called a point of order 9  a \emph{coincidence point}). In fact, 
It follows from the definition of the group law that we have
$2p\oplus q= 2q\oplus r = 2r\oplus s = 0.$ This immediately implies that $9p = 0$ if and only if $p = s$ (this explains why the classics called a point of order 9  a \emph{coincidence point}). The triangle formed by the lines $\overline{p,q}, \overline{q,r}, \overline{r,p}$  is inscribed and circumscribed to $E$. Following Halphen \cite{H}, we will use this observation to find the locus of points of order 9. 

The tangent line of $E$ at  $p = (x_0,y_0,z_0)$  has the equation
$$(x_0^2+ty_0z_0)x+(y_0^2+t x_0z_0)y+(z_0^2+t x_0y_0)z = 0,$$
where we assume that $E = E_{3t}$.  The point  $q = (x_0,\e y_0,\e^2 z_0)$ lies on $E$  because $(x_0,y_0,z_0)\in E$, it also lies on the tangent line at $p$ if  $p = (x_0,y_0,z_0)$ satisfies the equation
\begin{equation}\label{B1}
B_1: x^3+\e y^3+\e^2 z^3 = 0.
\end{equation}
If $p$  satisfies this equation, then $q$ also satisfies this equation,  hence 
$r = (x_0,\e^2 y_0, \e z_0)$ lies on the tangent at $q$ and again satisfies \eqref{B1}. Next time we repeat this procedure we get back to the original point $p$. Hence we see that any point in $B_1\cap E$ is a point of order 9.  Now we apply the elements of the Hessian group to the curve $B_1$ to get the remaining cubic curves $B_2,\ldots,B_8$. Notice that the stabilizer of $B_1$ in the Hessian group is generated by $\Gamma$ and $g_4$. It is a Sylow 3-subgroup of the Hessian group isomorphic to a semi-direct product $\Gamma\rtimes \bbZ/3\bbZ$.

To check the last assertion it is enough, using the $G_{216}$-action, to consider one of the curves $B_i$. For example, we see that the triangle $T_1$ of equation $xyz = 0$ is an inflection triangle of the curve $B_1$ and the triangles $T_2,T_3,T_4$ are inscribed and circumscribed to $B_1$. 
More precisely we have the following configuration:\\
\noindent i) $B_i$ and $B_{i+4}$ have $T_i$ as a common inflection triangle and they intersect in the 9 vertices of the other triangles;\\
ii) $B_i$ and $B_j$, $i\not=j$, $i,j\leq 4$,  intersect in the 3 vertices of  a triangle $T_k$ and they are tangent in the 3 vertices of $T_{\ell}$ with $k,\ell\not\in\{i,j\}$;\\
iii) $B_i$ and $B_{j+4}$, $i\not=j$, $i,j\leq 4$, intersect similarly with $k$ and $\ell$ interchanged.\\
For example, $B_1$ and $B_2$ intersect in the vertices of $T_3$ and are tangent in the vertices of $T_4$, while $B_1$ and $B_6$ intersect transversally on $T_4$ and are tangent on $T_3$.

\end{proof}

We will call the cubics $B_i$ the \emph{Halphen cubics}. Observe  that the element $g_0$ from the Hessian group sends $B_i$ to $B_{i+4}$. We  will call the pairs $(B_i,B_i' = B_{i+4})$ the \emph{pairs of Halphen cubics} and we will denote with $q_i,q'_i=q_{i+4}$ the corresponding pairs of points in $\Ps 2$.

It can be easily checked that the projective transformations $g_3,g_4$ act on the Halphen cubics as follows (with the obvious notation):
$$g_3: (121'2')(434'3'),\ g_4: (243)(2'4'3').$$

\begin{Rem}\label{gamma}
The linear representation of $\Gamma$ on the space of homogeneous cubic polynomials decomposes into the sum of one-dimensional eigensubspaces. The cubic polynomials defining $B_i$ together with the polynomials $xyz,\ x^3+y^3+z^3$ form a basis of  eigenvectors.
Moreover, note that the cubics $B_i$ are  equianharmonic cubics. In fact, they are all projectively equivalent to $B_1$, which is obviously isomorphic to the Fermat cubic. We  refer to \cite{Aure1}, \cite{Aure2} where the Halphen cubics play a role in the construction of bielliptic surfaces in $\bbP^4$. 
\end{Rem}

\begin{Rem} According to G. Halphen \cite{H} the rational map $\gamma\circ\beta^{-1}\circ p:\bbP^2 -\to \bbP^2$ can be explicitly given  by
$$(x,y,z) \mapsto (P_2'P_3'P_4', P_2P_3P_4,xyzP_1P_1'),$$
where $P_i, P'_i$ are the polynomials defining $B_i, B'_i$ as in Proposition \ref{hc}. His paper \cite{H}, besides many other interesting results, describes the locus of $m$-torsion points of nonsingular members of the Hesse pencil (see \cite{Frium} for a modern undertaking of this problem) %apparently without knowledge about Halphen's work). 
\end{Rem}

\begin{Rem} In characteristic 3 the cyclic group of projective transformations generated by $g_1$ acts on nonsingular members of the Hesse pencil as translation by 3-torsion points with the zero point taken to be $(1,-1,0)$. The polynomials 
$$(X,Y,Z,W) = (x^2y+y^2z+z^2x, xy^2+yz^2+zx^2, x^3+y^3+z^3,xyz)$$
are invariant with respect to $g_1$ and map $\bbP^2$ onto a cubic surface in $\bbP^3$ given by the equation (see \cite{Ran}, (3.1))
\begin{equation}\label{cubic}
X^3+Y^3 +Z^2W = XYZ.
\end{equation}
Among the singular points of the cubic surface,  $(0,0,0,1)$ is a rational double point of type $E_6^{(1)}$ in Artin's notation \cite{Artin}. The image of the member $E_\lambda$ of the Hesse pencil is the  plane section
$Z +\lambda W = 0$. Substituting in  equation \eqref{cubic}, we find that the image of this pencil of plane sections under the projection from the singular point is the  Hesse pencil. The parameter $\lambda$  of the original pencil and the new parameter $\lambda'$ are related by $\lambda = \lambda'{}^3.$
\end{Rem}
 
 \section{The 8-cuspidal sextic}\label{sextic}
Let $C_6$  be the sextic curve with equation $\Phi_6 = 0$, where $\Phi_6$ is the degree six invariant of the Hessian group. This is a smooth curve and it is immediately checked that it does not contain the vertices of the inflection triangles $T_1,\ldots, T_4$ given in \eqref{vertices} as well as  the base points of the Hesse pencil. 

This shows that the   preimage $\bar{C}_6=\pi^{-1}(C_6)$ of $C_6$ in the surface $S(3)$ is isomorphic to $C_6$ and the group $\Gamma$ acts on $\bar{C}_6$ freely. The orbit space $\bar{C}_6/\Gamma$ is a smooth curve of genus 2  on $S(3)/\Gamma$ which does not pass through the singular points and does not contain the orbit of the section $\pi^{-1}(p_0)$. Its   preimage under $\sigma$ is a smooth curve $\bar{C}_6'$ of genus 2  on $S(3)'$  that intersects a general fibre of the Hesse pencil at 2 points. Observe that the curve $C_6$ is tangent to each Halphen cubic $B_i,B_i'$ at a $\Gamma$-orbit of 9 points. In fact, it is enough to check that $C_6$ is tangent to one of them, say  $B_1$, at some point. We have
$$x^6+y^6+z^6-10(x^3y^3+x^3z^3+y^3z^3) = (x^3+y^3+z^3)^2-12(x^3y^3+x^3z^3+y^3z^3)$$
$$=  -3(x^3+y^3+z^3)^2+4(x^3+\e y^3+\e^2 z^3)(x^3+\e^2 y^3+\e z^3).$$
This shows that the curves $B_1$ and $B_1'$ are tangent to $C_6$ at the points where $C_6$ intersects the curve $E_0:x^3+y^3+z^3 = 0$.
The map $\pi':S(3)'\to \bbP^2$ blows down the curves $B_i,\ i = 1,\ldots,8,$ to the base points $q_1,\ldots,q_8,$ of the Hesse pencil.  Hence the image $C_6'$ of $\bar{C}_6'$ in $\bbP^2$ is a curve of degree 6 with cusps at the points $q_1,\ldots,q_8$.
 
\begin{Prop} The 8-cuspidal sextic $C_6'$  is projectively equivalent to the sextic curve defined by the polynomial
\small{\[\Phi_6'(x,y,z)= (x^3+y^3+z^3)^2-36y^3z^3+24(z^4y^2+z^2y^4)
-12(z^5y+zy^5)-12x^3(z^2y+zy^2).\]}
%The restriction of the map $\phi':S(3)'\to \bbP^1$ to its   preimage $\bar{C}_6'$ on $S(3)'$ 
%The projection of $C_6'$ from $p_0$ is a degree six cover of $\Ps 1$ ramified over the zeroes of the discriminant $\Delta$ of the Hesse pencil and over the zeroes of the polynomial $B(t_0,t_1)$ defining harmonic cubics in the Hesse pencil (see $(8)$).
\end{Prop}

\begin{proof} In an appropriate  coordinate system the points $q_i$ have the same coordinates as the $p_i$'s. By using the action of the group $\Gamma$, we may assume that the sextic has cusps at $p_1,\ldots,p_8$. Let $V$ be the vector space of homogeneous polynomials of degree 6 vanishing at $p_1,\ldots,p_8$ with multiplicity $\ge 2$. If $S$ is the blowing-up of $q_1,\dots,q_8$ and $K_S$ is its canonical bundle, then $\bbP(V)$ can be identified with the linear system  $\mid  -2K_S\mid $. It is known that the linear system $|-2K_S|$ is of dimension 3 (see \cite{Demazure}) and defines a regular map of degree 2 from $S$ to $\bbP^3$ with the image a singular quadric.

 A basis of  $V$ can be found by considering the product of six lines among the 12 inflection lines. In this way one finds the following sextic polynomials
\begin{equation}\label{newlines}
\small{\begin{array}{l}
A_1 = yz(x+\e y+z)(x+y+\e z)(x+\e^2y+z)(x+y+\e^2z),\\ 
A_2 = yz(x+\e y+\e^2 z)(x+\e^2y+\e z)(x+y+\e z)(x+\e y+z),\\ 
A_3 = yz(x+\e^2y+z)(x+y+\e^2z)(x+\e y+\e^2 z)(x+\e^2y+\e z),\\ 
A_4 = (x+\e y+\e^2 z)(x+\e^2y+\e z)(x+y+\e z)(x+\e y+z)(x+\e^2y+z)(x+y+\e^2z).
\end{array}}
\end{equation}
A polynomial $P(x,y,z)$ defining the curve $C_6'$ is invariant with respect to the linear representation of the binary tetrahedral group $\bar{T} \cong \SL(2,\bbF_3)$ in $V$. 
This representation decomposes into the direct sum of the 3-dimensional representation isomorphic to the second symmetric power of the standard representation of $2.A_4$ in $\bbC^2$ and a one-dimensional representation spanned by $P(x,y,z)$.  
Applying $g_4$ we find that $$(A_1,A_2,A_3,A_4)\mapsto (\e^2A_2,\e^2A_3,\e^2A_1,A_4).$$
Thus $P(x,y,z) = \lambda (A_1+\e^2A_2+\e A_3)+\mu A_4$ for some constants $\lambda,\mu$. Now we apply $g_3$ and find $\lambda, \mu$ such that $P(x,y,z)$ is invariant. A simple computation gives the equation of $C_6'$. \end{proof}

%To check the second assertion, consider the projection of the curve $C_6'$ from the point $p_0=(0,1,-1)$ onto the harmonic polar  $L_0=\{y-z = 0\}$. This is a degree six cover branched over $C_6'\cap L_0$ and the $4$ points of intersection of $L_0$ with the inflection lines through $p_0$. 
%The points in $C_6'\cap L_0$ are of type $(1,y,y)$ where $1-8y^6-20y^3 = 0$ (note that these are the fixed points of the involution $g_2$ acting on $C_6'$).
%The intersections of $L_0$ with the inflection lines through $p_0$ are given by $(-2\e^i,1,1)$, $i=0,1,2$.
%The fixed points of $g_2$ on $C_6'$ are the points $(1,y,y)$ given by the equation
%$1-8y^6-20y^3 = 0$. 
%The map $\phi'$ is given by $(x,y,z)\mapsto (t_0,t_1) = (-xyz,x^3+y^3+z^3)$. 
%Let $t_1/t_0 =6u$,  then eliminating  $u$ from the equation $u = -(1+2y^3)/y^2$, we find $1-20u^3-8u^6= 0$. 
%Taking the isomorphism $L_0\rightarrow \Ps 1$, $(1,y,y)\mapsto (1,y)$ and homogenizing, we obtain the result.
%the polynomial $B$ from $(8)$ whose locus of zeroes is the set of parameters defining the six harmonic cubics in the Hesse pencil. 

\begin{Rem}
The geometry of the surface $S$, the blow-up of $\Ps 2$ at  $q_1,\ldots,q_8$, is well-known. We now present several birational models of this surface and  relations between them. 

The surface $S$ is a Del Pezzo surface of degree $1$ and admits a birational morphism $\psi:S\to \bar{S}$ onto a surface  in the weighted projective space $\bbP(1,1,2,3)$ given by an equation
\begin{equation}\label{relation}
-u_3^2+u_2^3+A(u_0,u_1)u_2+B(u_0,u_1)= 0,
\end{equation}
where $(u_0,u_1,u_2,u_3)$ have weights $1,1,2,3$ (see \cite{Demazure}). 
%The surface is the anti-canonical model of a Del Pezzo surface of degree 1.  
The morphism $\psi$ is an isomorphism outside of the union of the 8 lines $\ell_1,\ldots,\ell_8$ which correspond to  factors of  the polynomials $A_1,\ldots,A_4$ from  \eqref{newlines}. In fact, the map $\psi$ is a resolution of indeterminacy points of the rational map $\bar\psi:\bbP^2\dashrightarrow \bbP(1,1,2,3)$. It is  given by the formulas
$$(x,y,z) \mapsto (u_0,u_1,u_2,u_3) = (-xyz, x^3+y^3+z^3, \Phi_6'(x,y,z),P_9(x,y,z)),$$
where $P_9(x,y,z)= 0$ is the union of the line $\ell_0:y-z = 0$ and the 8 lines $\ell_1,\dots,\ell_8$. Explicitly,
$$P_9(x,y,z) =  yz(y-z)(x^6+x^3(2y^3-3y^2z-3yz^2+2z^3)+(y^3-yz+z^2)^3.$$
Up to some constant factors, the polynomials $A,B$ are the same as in \eqref{binary}.
The 8 lines are blown down to singular points of the surface. %Essentially the same equation of $\bar{S}$ was given in \cite{vG}.

The composition of $\bar\psi$ with the projection $(u_0,u_1,u_2,u_3)\mapsto (u_0^2,u_0u_1,u_1^2,u_2)$ gives the rational map $\bbP^2 \dashrightarrow \bbP^3$ defined by
$$(x,y,z)\mapsto (u_0,u_1,u_2,u_3) = (x^2y^2z^2,xyz(x^3+y^3+z^3), (x^3+y^3+z^3)^2,\Phi_6').$$
This is a $2.A_4$-equivariant map of degree 2 onto the quadric cone $u_0u_2-u_1^2 = 0$. The ramification curve is the line $y-z = 0$ and the branch curve is the intersection of the quadric cone and a cubic surface. This is a curve $W$ of degree 6 with 4 ordinary cuspidal singularities lying on the hyperplane $u_3 = 0$. 

\label{dp} Consider the rational map  $\phi = \gamma\circ \beta^{-1}\circ p:\bbP^2-\to \bbP^2$ from the diagram \eqref{diagram}. It follows from the description of the maps in the diagram that the   preimage of the Hesse pencil is a Hesse pencil, the  preimage of the curve $C_6'$ is the curve $C_6$, and the   preimage of the union of the lines $\ell_0, \ell_1,\ldots,\ell_8$ is the union of harmonic polars. This shows that the composition $\psi\circ \phi:\bbP^2- \to \bbP(1,1,2,3)$ can be given by the  formulas
$$(x,y,z) \mapsto (xyz, x^3+y^3+z^3, \Phi_6(x,y,z), \Phi_9(x,y,z)),$$
where $\Phi_9(x,y,z)$ is the invariant of degree 9 for the group $\overline{G}_{216}$ given in \eqref{invariants}. This agrees with a remark of  van Geemen in \cite{vG} that the polynomials 
$xyz, x^3+y^3+z^3, \Phi_6(x,y,z),$ and $\Phi_9(x,y,z)$ satisfy the same relation \eqref{relation} as the polynomials
$xyz, x^3+y^3+z^3,\Phi_6'(x,y,z),$ and $P_9(x,y,z)$. 
Using the standard techniques of invariant theory of  finite groups one can show that the polynomials $xyz, x^3+y^3+z^3, \Phi_6(x,y,z),$ and $\Phi_9(x,y,z)$ generate the algebra of invariants of  the Heisenberg group $\calH_3(3)$, the   preimage of $\Gamma$ in $\SL(3,\bbF_3)$. 
 With respect to different sets of generators the equations of $S$ were given in \cite{Benvenuti} and \cite{Verlinde}.
%In the set of generators
%$$ X = xyz, \ Y = x^3+y^3+z^3, \ Z = (x^3+\e y^3+\e^2 z^3)(x^3+\e^2 y^3+\e z^3),$$
%$$ W = (x^3+\e y^3+\e^2 z^3)^3$$
%the relation is
%\begin{equation}\label{verlinde}
%W^2+Y^3-27WX^3+WZ^3-3WYZ = 0
%\end{equation}
%(see \cite{Verlinde}).
%In a different set of generators
%$$X = 3xyz,\  Y = x^3+y^3+z^3,\  Z = x^6+y^6+z^6,\ W = x^3y^6+x^6z^3+y^3z^6,$$
%the relation is 
%\begin{equation}\label{verlinde2}
%8X^6+X^3(-48Y^3+72YZ+72W)+81\bigl((Y^2-Z)^3-4Y(Y^2-Z)W+8W^2\bigr) = 0
%\end{equation}
%(see  \cite{Benvenuti}). 
\end{Rem}
%The affine hypersurface given by equations \eqref{verlinde} or \eqref{verlinde2} is an affine orbifold Calabi-Yau threefold. It is isomorphic to the quotient of $\bbC^3$ by the Heisenberg group $\calH_3 = 3.3^2$ which acts symplectically on $\bbC^3$. This is a special case of the following construction. Let $X$ be a Gorenstein projective algebraic variety with big and nef anti-canonical sheaf $\omega_X^{-1}$ and finitely generated  anticanonical ring 
%$$A_X = \bigoplus_{n=0}^\infty H^0(X,\omega_X^{-n}).$$
%Assume also that the projective spectrum of $A_X$ has only quotient singularities. Then the affine spectrum of $A_X$ is a non-compact orbifold Calabi-Yau manifold. Examples of varieties $X$ with these properties  are weak Del Pezzo surfaces,  toric varieties defined by simplicial polytopes, projective spectra of algebras of invariants of a finite subgroup $G$ of $\SL(N)$ which contains the subgroup generated by the scalar matrix $e^{2\pi i/N}I_N$.  The Hilbert functions
%$$f(t) = \sum_{n=0}^\infty \dim  H^0(X,\omega_X^{-n})t^n$$
%occur in physics as generating functions for counting  BPS-states in superconformal quiver gauge theories in the string $AdS_5\times X^5$-theories (see \cite{Benvenuti},\cite{Verlinde}).
%\end{Rem}

Finally, we explain the geometric meaning of the intersection points of the  sextic curve $C_6$ with a nonsingular member $E_\lambda$  of the Hesse pencil. This set of intersection points is invariant with respect to the translation group $\Gamma$ and the involution $g_0$, thus its image on $C_6'= C_6/\Gamma$ consists of two points on the curve $E_\lambda$. These points lie on the line through the point $p_0$ because they differ by the negation involution $g_0$ on $E_\lambda$ in the group law with the zero point  $p_0$.

\begin{Prop} The curves $C_6'$ and $E_\lambda$ intersect at two points $p,q$ outside the base points  $p_1,\ldots,p_8$. These points lie on a line through $p_0$ which is the tangent line to the Hessian cubic $\He(E_\lambda)$ at $p_0$. The 18 points in $C_6\cap E_\lambda$ are the union of the two $\Gamma$-orbits of $p$ and $q$.
\end{Prop} 
\begin{proof} This is checked by a straightforward computation. By using MAPLE we find that the curves $C_6'$, $E_\lambda$ and the tangent line to $E_{\he(\lambda)}$ at $p_0$ have two intersection points.
\end{proof}

\section{A $K3$ surface with an action of $G_{216}$}
In the previous sections we introduced two plane sextics, $C_6$ and $C'_6$ which are naturally related to the Hesse configuration.
The double cover of $\bbP^2$ branched along any of these curves is known to be birationally isomorphic to a \emph{K3 surface}, i.e. a simply connected compact complex surface with trivial canonical bundle. In the following sections we will study the geometry of the $K3$ surfaces associated to $C_6$ and $C'_6$, in particular we will show how the symmetries of the Hesse configuration can be lifted to the two surfaces.  %We will denote  by $X$ and $ X'$ the $K3$ surfaces which are birational to the double covers of $\bbP^2$ branched along $C_6$ and $C_6'$ respectively. 
We start with presenting some basic properties of $K3$ surfaces and their automorphisms (see for example \cite{BPV}).\\

Since the canonical bundle is trivial, the vector space $\Omega^2(X)$ of holomorphic 2-forms on a $K3$ surface $X$ is one-dimensional.  
Moreover, the cohomology group $L = H^2(X,\bbZ)$ is known to be a free abelian group of rank 22. The cup-product equips $L$ with a structure of a \emph{quadratic lattice}, i.e. a free abelian group together with an integral quadratic form.  The quadratic form is unimodular and its signature is equal to  $(3,19)$. The sublattice $S_X\subset L$ generated by the fundamental cocycles of algebraic curves on $X$ is called the \emph{Picard lattice} and its orthogonal complement $T_X$ in $L$ is  the \emph{transcendental lattice} of $X$.

Any automorphism $g$ on $X$ clearly acts on $\Omega^2(X)$ and also induces an isometry $g^*$ on $L$ which preserves $S_X$ and $T_X$. 
An automorphism  $g$ that acts identically on $\Omega^2(X)$ is called \emph{symplectic}.  We recall here a result given in \cite{N3}.
\begin{Thm}\label{auto}
Let $g$ be an automorphism of finite order on a $K3$ surface $X$.
\begin{itemize}
\item[i)] If $g$ is symplectic then $g^*$ acts  trivially on $T_X$
and its fixed locus is a finite union of points.
%Moreover $o(g)\leq 8$ and its fixed locus is the union of $n$ points, where $n$ depends on the order of $g$ as follows:
%$$\begin{array}{c|c|c|c|c|c|c|c|c}
%o(g) & 1 & 2 & 3 & 4 & 5 & 6 & 7 & 8\\
%\hline
%n& 24 & 8 & 6 & 4 & 4 & 2 & 3 & 2  
%\end{array}\ ;$$ 
\item[ii)]
If $g$ is not symplectic then $g^*$ acts on $\Omega^2(X)$ as the multiplication by a primitive $r$-th root of unity and its eigenvalues on  $T_X$ are the primitive $r$-th roots of unity. \end{itemize}
\end{Thm}
%If this happens then the set of fixed points of any non-trivial element of the group is finite and its cardinality  depends only on the order of such element (see \cite{N3}).
%In the following we will denote  by $X$ and $ X'$ the $K3$ surfaces which are birational to the double covers of $\bbP^2$ branched along $C_6$ and $C_6'$ respectively. 
Let $q:X\lra \Ps 2$ be the double cover branched along $C_6$.
We now prove that the action of the Hessian group on the projective plane lifts to an action on $X$.
We denote by $Q_8$ be the 2-Sylow subgroup of $\SL(2,\bbF_3)$, isomorphic to the quaternion group.
\begin{Prop}\label{act}
The Hessian group $G_{216}$ is isomorphic to a group of automorphisms of the $K3$ surface $X$.  Under this isomorphism, any automorphism in the normal subgroup $H_{72} = \Gamma\rtimes Q_8$ is symplectic.
\end{Prop}
\begin{proof} 
The double cover $q:X\to \bbP^2$ branched along the curve $C_6$ can be defined by the equation 
$$w^2+\Phi_6(x,y,z) = 0,$$
considered as a weighted homogeneous polynomial with weights $(1,1,1,3)$. Thus we can consider $X$ as a hypersurface of degree 6 in the weighted projective space $\bbP(1,1,1,3)$.

Let $G_{216}'$ be the preimage of $G_{216}$ in $\SL(3,\bbC)$ considered  in section \ref{hessegroup} and $g_i', i = 1,\ldots,4,$ be the lifts of the generators $g_i$ in $G_{216}'$.  It is checked immediately that  the generators $g_1',g_2',g_3'$ leave the polynomial $\Phi_6$ invariant and $g_4'$ multiplies $\Phi_6$ by $\e^2$.  Thus the group $G_{216}'$ acts on $X$ by the formula
$$g_i(x,y,z,w) = (g_i'(x,y,z),w), \ i \ne 4, \quad g_4(x,y,z,w) = (g_4'(x,y,z),\e w).$$
The kernel of $G_{216}'\to G_{216}$ is generated by the scalar matrix $(\e,\e,\e)$, which acts as the identity transformation on $X$.  Then it is clear that the induced action of $G_{216}$ on $X$ is faithful. 

The subgroup $H_{72}$ of $G_{216}$ is generated by the transformations
$g_1,g_2,g_3, g_4g_3g_4^{-1}$. To check that it acts symplectically on $X$ we recall that the space of holomorphic 2-forms on a hypersurface $F(x_0,\ldots,x_n)$ of degree $d$ in $\bbP^n$ is generated by the residues of the meromorphic $n$-forms on $\bbP^n$ of the form
$$\omega = \frac{P}{F}\sum_{i=0}^n(-1)^ix_idx_1\wedge\cdots \wedge \widehat{dx_i}\wedge\cdots \wedge dx_n,$$
where $P$ is a homogeneous polynomial of degree $d-n-1$. This is easily generalized to the case of hypersurfaces in a weighted projective space $\bbP(q_0,\ldots,q_n)$. In this case the generating forms are
$$\omega = \frac{P}{F}\sum_{i=0}^n(-1)^iq_ix_idx_1\wedge\cdots \wedge \widehat{dx_i}\wedge\cdots \wedge dx_n,$$
where $\deg P = d-q_0-\cdots-q_n$. In our case $d=q_0+q_1+q_2+q_3=6$, hence there is only one form, up to proportionality. It is given by
$$\omega = \frac{xdy\wedge dz\wedge dw-ydx\wedge dz\wedge dw+zdx\wedge dy\wedge dw-3wdx\wedge dy\wedge dz}{w^2+\Phi_6(x,y,z)}.$$
It is straightforward to check that the generators of $H_{72}$ leave this form invariant (cf. \cite{M}, p. 193).

%However, $g_3^2$ acts by $(x,y,z,w) \mapsto (-x,-z,-y,w) = (x,z,y,-w)$. It is immediately seen that it has only finitely many fixed points on $X$ (in fact, 8 of them). This shows that $g_3$ has only isolated fixed points on $X$ and hence acts symplectically. Thus $\chi(g_3) = 1$, hence $\chi$ is trivial and $H_{72}$ acts symplectically.  
\end{proof}

\begin{Rem}
The action of $\Gamma\rtimes Q_8$ appears as Example 0.4  in the paper of S. Mukai \cite{M} containing the classification of maximal finite  groups of symplectic automorphisms of complex $K3$ surfaces.
\end{Rem}
%As remarked in \cite{OZ} all $K3$ surfaces with maximum symplectic action have maximum Picard number.
%\begin{Lem}
%The $K3$ surface $X$ is singular and the transcendental lattice is a free $\Z[\e]$-module.
%\end{Lem}
%\proof
Let $P_i,P'_i$, $i=1,\dots,4,$ be the polynomials defining the cubics $B_i, B_{i+4}$ as given in section $5$ and $F_i$ be the equations of the equianharmonic cubics in the Hesse pencil:
%$$F_1(x,y,z)=x^3+y^3+z^3,$$
$$ F_i(x,y,z)=x^3+y^3+z^3+\alpha_i xyz,\ \ \ i=1,\dots,4,$$
where $\alpha_1=0$ and $\alpha_i=6\e^{2-i}$ for $i=2,3,4$ (see section 2).
 
\begin{Prop}\label{mod}
The K3 surface $X$ is isomorphic to the hypersurface of bidegree $(2,3)$ in  $\Ps 1\times \Ps 2$ with equation
\begin{equation}\label{new1}
u^2 P_i(x,y,z)+ v^2 P'_i(x,y,z)+\sqrt{3}uv  F_i(x,y,z)=0
\end{equation}
for any $i=1,\dots,4$.
\end{Prop}
\begin{proof}
As noticed in the previous section we can write
%\label{det}
$$
 \Phi_6=\det\left(\begin{array}{cc}
2P_1& \sqrt{3} F_1\\
\sqrt{3}F_1 & 2P'_1\end{array}\right
)=-3F_1^2+4P_1P'_1.
$$ 
The $K3$ surface $Y$ given by the bihomogeneous equation  of bidegree $(2,3)$ in $\Ps 1\times \Ps 2$
\begin{equation}\label{new}
u^2 P_1(x,y,z)+ v^2 P'_1(x,y,z)+\sqrt{3}uv  F_1(x,y,z)=0
\end{equation}
is a double cover of $\bbP^2$ with respect to the projection to the second factor and its branch curve is defined by $\Phi_6=0$. Thus $Y$ is isomorphic to the $K3$ surface $X$. By acting on equation (\ref{new}) with the Hessian group $G_{216}$ we find analogous equations for $X$ in $\Ps 1\times \Ps 2$ in terms of the polynomials $P_i, P'_i$ and $F_{i}$ for $i=2,3,4$.
%In other words, $X$ has four  different models as a divisor in  $\Ps 1\times \Ps 2$.
\end{proof}
An important tool for understanding the geometry of $K3$ surfaces is
the study of their elliptic fibrations.
%We recall that the \emph{multisection index} of an elliptic fibration  is the minimum positive intersection number between a fibre and a divisor. 
We recall that  the fibration is called \emph{jacobian} if it has a section.
%Its elements can be represented by sections and the translation automorphism of the general fibre corresponding to a $K$-rational point can be extended to an automorphism of the surface.

%\noindent Proposition \ref{mod} shows the existence of some natural elliptic fibrations on $X$.
\begin{Prop}\label{fib}
The $K3$ surface $X$ has $4$ pairs of elliptic fibrations   
$$\xymatrix{\frakh_i,\frakh_i':X \ar[r]& \bbP^1,\ \ \ i=1,\dots,4}$$ 
with the following properties:
\begin{itemize}
\item[$a)$] $\frakh_i$ and  $\frakh_i'$ are exchanged by the covering involution of $q$ and $G_{216}$ acts transitively on $\frakh_1,\dots, \frakh_4$; 
\item[$b)$] the $j$-invariant of any smooth fibre of $\frakh_i$ or $\frakh_i'$ is equal to zero;
\item[$c)$]  each fibration has $6$ reducible fibres of Kodaira's type $IV$, i.e. the union of three smooth rational curves intersecting at one point. The singular points in the reducible fibres of $\frakh_i$ and $\frakh_i'$ are mapped by $q$ to the vertices of the triangle $T_i$;
\item[$d)$]  each fibration is jacobian.
\end{itemize}
\end{Prop}

\proof  Consider the equations \eqref{new1} for $X$ in $ \Ps 1\times \Ps 2$. The projections on the first factor $\frakh_i:X\to \Ps 1$,  $i=1,\dots,4$ are elliptic fibrations on $X$ since the fibre over a generic point $(u,v)$ is a smooth plane cubic.
A second set of elliptic fibrations on $X$ is given by 
 $\frakh_i'=\frakh_i\circ \sigma$, where $\sigma$ is the covering involution of $q$. 
%Note that the plane cubic $B_i$ is the projection to $\Ps 2$ 
 %intersects it at 9 points and clearly projects to the same cubic on the plane. 
%It is everywhere tangent to $C_6$ and splits into  $F\cup \sigma(F)$ under the cover $X\to \bbP^2$. 
%By considering the other equations  in (\ref{new1}) we find $3$ other pairs of elliptic fibrations on $X$: $\frakh_i$ and $\frakh'_i=\frakh_i\circ \sigma$, $i=2,3,4$. This gives $a)$.
Since all these fibrations are equivalent modulo the group generated by  $\sigma$ and $G_{216}$, it will be enough to prove properties $b)$, $c)$ and $d)$ for $\frakh_1$.  

The fibre of $\frakh_1$ over a point $(u,v)$ is isomorphic to the plane cubic defined by equation \eqref{new}. This equation can be also written in the form
\begin{equation}\label{neweq}
(u^2+v^2+\sqrt{3}uv)x^3+(\e u^2+\e^2v^2+\sqrt 3uv)y^3+(\e^2u^2+\e v^2+\sqrt{3}uv)z^3=0.\end{equation}
Hence it is clear that all smooth fibres of $\frakh_1$ are isomorphic to a Fermat cubic i.e. they are equianharmonic cubics. This system of plane cubics contains exactly $6$ singular members corresponding to the vanishing of the coefficients at  $x^3$, $y^3$ and $z^3$.  Each of them  is equal to  the union of three  lines meeting at one point and define six singular  fibres of type $IV$ of the elliptic fibration $\frakh_1$.  The singular  points of  these reducible fibres are the inverse images  of the vertices $v_0$, $v_1$, $v_2$ of the triangle $T_1$ under the map $q$ (see section 2, (\ref{vertices})).
This proves assertions $b)$ and $c)$.

It remains to show that  the elliptic fibration $\frakh_1$ has a section. We thank N. Elkies for finding explicitly such a section. It is given by 
\begin{equation}\label{elkies}
(x,y,z)=((1-\epsilon)u-d_0v,(1-\epsilon)u-d_1v, (1-\epsilon)u-d_2v),
%(x,y,z) = ((1-\epsilon) u+\sqrt[3]{4}i v, (1-\epsilon)u-(1+\epsilon)iv, (1-\epsilon)u-\epsilon\sqrt[3]{2}iv).
\end{equation}
where $(d_0,d_1,d_2)=ig_3(-\sqrt[3]{2}, \epsilon+1, \epsilon \sqrt[3]{2})$.
 \qed
\begin{Rem}\label{rem4.2}
 Consider a map  
$$\xymatrix{\bbP^2\ar[r]& \bbP^2},\ (x,y,z)\longmapsto (x^3,y^3,z^3).$$
The image of  the curve $C_6$ is a conic $K$ and the preimages of the tangent lines to $K$  are plane cubics that everywhere tangent to $C_6$.
The map induces a degree $9$ morphism from $X$ to $\bbP^1\times \bbP^1$ isomorphic to the double cover of $\Ps 2$ branched along $K$.  The projections to the two factors give the fibrations $\frakh_1$ and $\frakh_1'$.
 
 Note that each family of everywhere tangent cubics to $C_6$ corresponds to an even theta characteristic $\theta$ on $C_6$ with $h^0(\theta)=2$.
\end{Rem}

Let $\pi:Y\to \bbP^1$ be a jacobian elliptic fibration on a K3 surface $Y$. The  fibre of $\pi$ over the generic point $\eta$ is an elliptic curve $Y_\eta$ over the field of rational functions $K$ of $\bbP^1$. The choice of a section $E$ of $\pi$ fixes a $K$-rational point on $Y_\eta$ and hence allows one to find a birational model of $Y_\eta$ given by a Weierstrass equation $y^2-x^3-ax-b = 0$, where 
$a,b\in K$. The construction of the Weierstrass model can be ``globalized'' to obtain the following birational model of $Y$ (see \cite{CD}).

\begin{Prop} There exists a birational morphism $f:Y\to W$, where $W$ is a hypersurface in the weighted projective space $\bbP(1,1,4,6)$ given by an equation of degree 12
$$y^2-x^3-A(u,v)x-B(u,v) = 0,$$
where $A(u,v), B(u,v)$ are binary forms of degrees $8$ and $12$ respectively. Moreover
\begin{enumerate}[$\bullet$]
\item The image of the section $E$ is the point $p = (0,0,1,1)\in W$. The projection $(u,v,x,y)\mapsto (u,v)$ from $p$ gives an elliptic fibration $\pi':W'\to \bbP^1$ on the blow-up $W'$ of $W$ with center at $p$.  
It  has a section defined by the exceptional curve $E'$ of the blow-up. 
\item The map $f$ extends to a birational morphism $f':Y\to W'$ over $\bbP^1$ which maps $E$ onto $E'$ and blows down irreducible components of fibres of $\pi$ which are disjoint from $E$ to singular points of $W'$. 
\item Each singular point of  $W'$ is a double rational point of type $A_n, D_n,E_6,E_7$ or $E_8$.  A singular point  of type $A_n$ corresponds to a fibre of $\pi$ of Kodaira type $I_{n+1}$, $III$ (if $n= 1$), or $IV$ (if $n= 2$). A singular point of type $D_n$  corresponds to a fibre of type $I_{n+4}^*$. A singular point of type  $E_6, E_7, E_8$ corresponds to a fibre of type $IV^*, III^*, II^*$ respectively.
\end{enumerate}
\end{Prop}

The elliptic surface $W$ is determined uniquely, up to isomorphism,  by the elliptic fibration on $Y$.  It is called the \emph{Weierstrass model} of the elliptic fibration $\pi$.

It is easy to find the  Weierstrass model of the elliptic fibration $\frakh_1:X\to \bbP^1$ on our surface $X$. 

\begin{Lem} The Weierstrass model of the elliptic fibration $\frakh_1$ is given by the equation
$$y^2-x^3-(u^6+v^6)^2= 0.$$
\end{Lem}
\proof We know from Proposition \ref{fib} that the $j$-invariant of a general fibre of $\frakh_1$ is equal to zero. This implies that the coefficient $A(u,v)$ in the Weierstrass equation is equal to zero. We also know that the fibration has 6 singular fibres of type $IV$ over the zeroes of the polynomial
$$(u^2+v^2+\sqrt{3}uv)(\e u^2+\e^2v^2+\sqrt 3uv)(\e^2u^2+\e v^2+\sqrt{3}uv) = u^6+v^6.$$
Since each of the fibres is of Kodaira type $IV$, the singularity of $W$ over a root of $u^6+v^6$ must be a rational double point of type $A_3$, locally isomorphic to the singularity $y^2+x^3+z^2$. This easily implies that the binary form $B(u,v)$ is equal to $(u^6+v^6)^2$ up to a scalar factor which does not affect the isomorphism class of the surface. 
\qed

\begin{Lem}\label{o3}
Let $Y$ be a K3 surface with Picard number $20$ having a non-symplectic automorphism of order $3$. Then the intersection matrix of $T_Y$ with respect a suitable basis is given by  
\begin{equation}\label{matrix}
A_2(-m)=\left(
\begin{array}{cc}
2m& m\\
m& 2m\end{array} \right 
), \end{equation}
for some $m\in \Z$, $m>0$.  
\end{Lem}
\proof
Let $f$ be a non-symplectic automorphism of order $3$ on $Y$. By Theorem \ref{auto}, ii),  $f^*$ acts on $T_Y$ as an isometry with eigenvalues $\e, \e^2$.
Let $x\in T_Y$, $x\not=0$,  then 
$$0=(x+f^*(x)+(f^*)^2(x),f^*(x))=2(x,f^*(x))+x^2.$$ 
Note that $x^2=2m$ for some positive integer $m$ because the lattice $T_Y$ is even and positive definite.
Then the intersection matrix of $T_Y$ with respect to the basis  $x, -f^*(x)$ is $A_2(-m)$. See also Lemma 2.8 in \cite{OZ}.
\qed\\

%The automorphism $g_4(x,y,z,w)=(x,\e y,\e z,\e w)$ on $X=\{w^2+\Phi_6(x,y,z)=0\}$ clearly fixes the curve $\{x=0\}$ pointwisely. Hence, by Theorem \ref{auto}, $g_4$ has a non-symplectic action on $X$ and  $g^*_4$ acts on $T_X$ as an isometry with eigenvalues $\e, \e^2$. 
%Note that the lattice corresponding to $m = 1$ is isomorphic to the lattice $A_2(-1)$, where $A_2$ is the root lattice of type $A_2$ with the quadratic form multiplies by $-1$ (this sign confusing notation is traditionally used in algebraic geometry). Thus we see that $T_X \cong  A_2(-m)$.

\noindent The proof of the following theorem follows a suggestion of M. Schuett.

\begin{Thm}\label{k1}
The intersection matrix of the transcendental lattice of the $K3$ surface $X$  with respect to a suitable basis is:
$$A_2(-6)=\left(
\begin{array}{cc}
12&6\\
6& 12\end{array} \right 
). $$
 \end{Thm}
 
\proof

Consider the automorphism $\sigma$ of order $6$ of $X$ that acts on the Weierstrass model by the formula $(u,v,x,y) \mapsto (\eta u, v, \eta^2 x, \eta^3y)$, where $\eta = e^{\pi i/3}$. It is easy to see that $\sigma$ acts freely outside the union of the two nonsingular fibres $F_0, F_\infty$ over the points $(u,v) = (1,0)$ and $(0,1)$. The action of the cyclic group $G = \langle\sigma\rangle$ on each of the fibres is an automorphism of order $6$ such that $G$ has one fixed point, $ \langle\sigma^3\rangle$ has $4$ fixed points and $ \langle\sigma^2\rangle$ has $3$ fixed points. 

Let $X/G$ be  the orbit space. The images $\bar{F}_0$ and $\bar{F}_\infty$  of $F_0$ and $F_\infty$ in $X/G$ are smooth rational curves and $X/G$ has 3 singular points on each of these curves, of types $A_5, A_3$ and $A_2$. A minimal resolution of $X/G$ is a K3 surface $Y$. The elliptic fibration $\frakh_1$  on $X$ defines an elliptic fibration $\rho:Y\to \bbP^1$ with two fibres of type $II^*$, equal to preimages of $\bar{F}_0$ and $\bar{F}_\infty$ on $Y$, and one fibre of type $IV$, the orbit of the six singular fibres of $\frakh_1$. 

 It is easy to compute the Picard lattice $S_Y$ of $Y$. Its sublattice generated by irreducible components of fibres and a section of $\rho$ is isomorphic to $U\oplus E_8\oplus E_8\oplus A_2$, where $U$ is generated by a general fibre and a section. It follows from the Shioda-Tate formula  in \cite{Shioda2} that this sublattice coincides with $S_Y$ and  the  discriminant of its quadratic form is equal to $-3$. Since the transcendental lattice $T_Y$ is equal to the orthogonal complement of $S_Y$ in the unimodular lattice $L = H^2(X,\bbZ)$, this easily implies that $T_Y$ is a rank 2 positive definite even lattice with discriminant equal to $3$. There is only one isomorphism class of such a lattice and it is given by  $A_2(-1)$. 
 
 The transcendental lattices of the surfaces $X$ and $Y$ are related in the following way. By Proposition 5, \cite{Shioda3}, there is an isomorphism of abelian groups $(T_Y)\otimes \bbQ\cong (T_X)^G\otimes \bbQ$, defined by taking the inverse transform of transcendental cycles under the rational map $X\to Y$.  Since $G$ acts symplectically on $X$, we have $(T_X)^G = T_X$.   Under this map the intersection form is multiplied by the degree of the map equal to $6$. This implies that $T_X$ has rank two and contains $T_Y(6)\cong A_2(-6)$ as a sublattice of finite index.  

Note that the automorphism $g_4(x,y,z,w)=(x,\e y,\e z,\e w)$ of $X $ clearly fixes the curve $\{x=0\}$ pointwisely. Hence  $g_4$ is non-symplectic by Theorem \ref{auto}. It follows from Lemma \ref{o3} and the previous remarks that  $T_X\cong A_2(-m)$. Hence we only need to determine the integer $m$. 
As we saw in above, $T_X$ contains a sublattice isomorphic to $A_2(-6)$, hence $m\in\{1,2,3,6\}$.
%gives us the following possibilities: $T_X \cong A_2(-1), A_2(-2), A_2(-3)$ and $A_2(-6)$. 
We now exclude all possibilities except the last one. 
 
 The K3 surface with $T_Y \cong A_2(-1)$ was studied in \cite{N}, in particular all  jacobian elliptic fibrations on $Y$ are classified in Theorem 3.1. Since no one of these fibrations has the same configuration of singular fibres as $\frakh_1$ (see Proposition \ref{fib}), this excludes the case $m=1$. 
 
 The K3 surface with $T_Y\cong A_X(-2)$ is isomorphic to the Kummer surface from Theorem \ref{kummer} below. All its jacobian fibrations are described in \cite{N}, Theorem 3.1 (Table 1.1) and, as in the previous case, none of them has $6$ fibres of type $IV$. This excludes case $m=2$.
 
Finally, a direct computation shows that  $A_2(-3)$ does not contain a sublattice isomorphic to $A_2(-6)$. In fact, since the equation $x^2+y^2+xy=2$ has no integral solutions, then $A_2(-3)$
does not contain any element with self-intersection $12$.
%(it has a sublattice of index 2 but it is not isomorphic to $A_2(-3)$). We leave this computation to the reader. 
This completes the proof of our theorem.

\qed
 
\noindent We conclude this section giving another model for the surface $X$.
\begin{Prop}
The $K3$ surface $X$ is birational to the double cover of $\, \Ps 2$ branched along a sextic with $8$ nodes which admits a group of linear automorphisms isomorphic to $A_4$.
\end{Prop}
\proof
The lift of the involution $g_0$ to the cover $X=\{w^2+\Phi_6(x,y,z)=0\}$ given by
$$\tilde g_0(x,y,z,w)=(x,z,y,w)$$
is a non-symplectic involution.  The fixed locus of $\tilde g_0$ is the genus two curve $\tilde{L}_0$ which is the double cover of the harmonic polar $L_0=\{y-z=0\}$ branched along $L_0\cap C_6$. The quotient surface $R= X/(\tilde g_0)$  is a Del Pezzo surface  of degree 1, the double cover of $\bbP^2/(g_0) \cong Q$, where $Q$ is the quadratic cone with vertex equal to the orbit of the fixed point $p_0 = (0,1,-1)$ of $g_0$. 
%The pencil of lines through $p_0$ descends to the ruling of $Q$ and its   preimage on $R$ is the elliptic pencil $|-K_{R}|$ with base point at the   preimage of the vertex. Blowing-up the base point we obtain a rational elliptic surface $f: R'\to \bbP^1$ and the branch curve $B$ of $X\to R$ defines a  2-section  of $f$.  The fibration $f$ has 12 irreducible singular fibres corresponding to the lines through $p_0$ which are tangent to the curve $C_6$  at 2 distinct points forming a $g_0$-orbit.  There are also 6 tangents to $C_6$  at the points where $C_6$ intersects $L_0$, these define 6 members of the elliptic pencil which are isomorphic to equianharmonic cubics. 

Let $b:R\to \bbP^2$ be the blowing-down of $8$ disjoint $(-1)$-curves on $R$ to points $s_1,\ldots,s_8$ in $\bbP^2$. The pencil of cubic curves through the eight points is the image of the elliptic pencil $|-K_R|$ on $R$. Note that the stabilizer of the point $p_0$ in the Hessian group is isomorphic to $2.A_4$ with center equal to $(g_0)$, thus the group $A_4$ acts naturally on $R$ and  on the elliptic pencil $|-K_R|$. The curve $B\in |-2K_R|$ is an $A_4$-invariant member of the linear system $|-2K_R|$ and  $b(B)$ is a plane sextic with $8$ nodes at the points $s_1,\ldots,s_8$.  Thus we see that $X$ admits $9$ isomorphic  models as a double cover of the plane branched along a $8$-nodal sextic  with a linear action of $A_4$.
\qed

\begin{Rem} In \cite{KOZ} the authors study a $K3$ surface birationally isomorphic to the double cover of $\bbP^2$ branched along the union of two triangles from the Hesse pencil. This surface has transcendental lattice of rank $2$ with intersection matrix $\big(\begin{smallmatrix}6&0\\
0&6\end{smallmatrix}\big)$ and it admits a group of automorphisms isomorphic to $A_6\rtimes \bbZ/4\bbZ$. 
\end{Rem}

\section{A $K3$ surface with an action of $SL(2,\mathbb F_3)$}
We now study the $K3$ surface which is birational to the 
double cover of $\bbP^2$ branched along the sextic 
$C_6'$ defined by $\Phi_6'= 0$. 

We recall that $C_6'$ has $8$ cusps in the base points $q_1,\dots,q_8$ of the Hesse pencil. The double cover of $\Ps 2$ branched along $C_6'$  is locally isomorphic to $z^2+x^2+y^3 = 0$ over each cusp of $C_6'$, hence it has $8$ singular points of type $A_2$ (see \cite{BPV}).   
It is known that the minimal resolution of singularities of this surface is a $K3$ surface and that the exceptional curve over each singular point of type $A_2$ is the union of two rational curves intersecting in one point.

In this section we will study the properties of this $K3$ surface, which will be denoted by $X'$.  

\begin{Prop}\label{quo} The $K3$ surface $X'$ is birationally isomorphic to the quotient of the $K3$ surface $X$ by the subgroup $\Gamma$ of $G_{216}$. In particular, the group $SL(2,\mathbb F_3)$ is isomorphic to a group of automorphisms of $X'$.
\end{Prop}
\begin{proof}  
The minimal resolution of the double cover of $\Ps 2$ branched along $C_6'$ can be obtained by first resolving the singularities of $C_6'$ through the morphism $\gamma:S\to \bbP^2$ from  diagram \eqref{diagram} and then taking the double cover $\bar q':X'\to S$ branched over the proper transform $\bar{C}_6'$ of $C_6'$ (\cite{BPV}). 
Since $\beta^{-1}(p(C_6)) = \bar{C}'_6$ we have the following commutative diagram
\begin{equation}\label{diagram2}
\xymatrix{
X \ar[d]_{q} \ar[r] ^{\tilde r}&X/\Gamma\ar[d]&X'\ar[d]_{\bar q'}\ar[l]_{\tilde \sigma}\ar[rd]^{q'}& \\
  \Ps 2 \ar[r]^p& \Ps 2/\Gamma&\ar[l]_\beta S \ar[r]^{\gamma}&\Ps 2,}
 \end{equation}
where $q$ and $q'$ are the double covers branched along $C_6$ and $C_6'$ respectively,  $\tilde r$ and $p$ are  the natural quotient maps,  $\tilde \sigma$ is a minimal resolution of singularities and the bottom maps are as in diagram \eqref{diagram}. 
This gives the first statement. The second one follows from the isomorphism $G_{216}/\Gamma \cong SL(2,\mathbb F_3)$. %Note that the there are exactly $12$ points with non-trivial stabilizer subgroup in $\Gamma$: these are  the vertices of the triangles in the Hesse pencil (see ...).
%The vertices belong to $4$ orbits for the action of $\Gamma$ and the local action of the stabilizer at these points is of type $diag(\e,\e^2)$. Hence
%the quotient $X/\Gamma$ has $4$ singular points of type $A_2$ in the image of these orbits. 
\end{proof}  
%Note that the quotient group $M_9/\Gamma\cong Q_8$ acts symplectically on $X'$.
%A computation similar to that given in the proof of Proposition \ref{} gives that the invariant lattice has %rank $5$,
%In particular the Picard number of $X'$ is at least $18$ since the transcendental lattice has even rank. 
%In fact, the N\`eron Severi lattice contains the lattice $\langle 2\rangle\oplus 8A_2$.
\begin{Rem}\label{x/gamma}
The points in $X$ with non trivial stabilizer for the action of $\Gamma$ are exactly the $24$ preimages by $q$ of the vertices of the triangles in the Hesse pencil. 
In fact these points belong to $8$ orbits for the action of $\Gamma$ and give $8$ singular points of type $A_2$ in the quotient surface $X/\Gamma$ (see Proposition \ref{quo2}).
\end{Rem}
\noindent We now describe some natural elliptic fibrations on the surface $X'$.

\begin{Prop}\label{iac}
The pencil of lines through each of the cusps of $C_6'$ induces a jacobian fibration on $X'$ with $3$ singular fibres of Kodaira's type $I_6$ and one of type $I_3$ (i.e. cycles of $6$ and $3$ rational curves respectively).
\end{Prop}
\proof
Let $p$ be a cusp of $C_6'$ and $h_p$ be the pencil of lines through $p$. The generic line in the pencil intersects $C_6'$ in $p$ and $4$ other distinct points, hence its preimage in $X'$ is an elliptic curve. Thus $h_p$ induces an elliptic fibration $\tilde h_p$ on $X'$.

The pencil $h_p$ contains $3$ lines through $3$ cusps and one line through $2$ cusps of $C_6'$ since the cusps of $C'_6$ are the base points of the Hesse pencil. 
The proper transform  of a line containing $3$ cusps is a disjoint union of two smooth rational curves. Together with the preimages of the cusps, the full preimage of such line in $X'$ gives a fibre of $\tilde h_p$ of  Kodaira's type $I_6$, described by the affine Dynkin diagram  $\tilde{A}_5$. 
Similarly, the preimage of a line containing 2 cusps gives a fibre of $\tilde h_p$ of  type $I_3$ (in this case the proper transform of the line does not split). 
Thus $\tilde h_p$ has three fibres of type $I_6$ and one of type $I_3$  .

The exceptional divisor over the cusp $p$ splits in two rational curves $e_1,e_2$ on $X'$ and each of them intersect each fibre of $\tilde h_p$ in one point i.e. it is a section of $\tilde h_p$.
\qed

\begin{Prop}\label{fib2}
The elliptic fibrations $\frakh_i,\frakh_i'$, $ i=1,\ldots,4$ on $X$ induce $8$ elliptic fibrations $\bar \frakh_i,\bar\frakh_i'$ on $X'$ such that
\begin{itemize}
\item[$a)$]  $\bar\frakh_i$ and  $\bar\frakh_i'$ are exchanged by the covering involution of $q'$ and $SL(2,\mathbb F_3)$ acts transitively on $\bar\frakh_1,\dots, \bar\frakh_4$;
\item[$b)$] the $j$-invariant of a smooth fibre of the elliptic fibration $\bar\frakh_i$ or $\bar\frakh_i'$ is equal to zero; 
\item[$c)$]  each fibration has two fibres of Kodaira's type $IV^*$ (i.e. $7$ rational curves in the configuration described by the affine Dynkin diagram $\tilde {E_6}$) and two of  type $IV$. 
 \end{itemize}
%The image of any of these fibrations by the map $q'$ is  a one-dimensional family of curves of degree $9$ in $\Ps 2$ with $8$ triple points in $q_1, \dots,q_8$ and $3$ cusps on $C_6'$.
 \end{Prop}
\proof
It will be enough to study  the fibration $\frakh_1$, since all  other fibrations are projectively equivalent to this one by the action of $G_{216}$ and $\sigma$.

Let $g_1, g_2$ be the generators of $\Gamma$ as in section 4.
The polynomials  $P_1, P'_1$ and $F_1$ are eigenvectors for the action of $\Gamma$ (Remark \ref{gamma}), hence it is clear from equation (\ref{new}) that $\Gamma$ preserves the elliptic fibration $\frakh_1$. In fact, $g_1$ acts on the basis of the fibration as an order three automorphism and fixes exactly the two fibres $E_1,E'_1$ such that $q(E_1)=B_1$ and $q(E'_1)=B'_1$.
The automorphism $g_2$ preserves each fibre of $\frakh_1$ and acts on it as an order $3$ automorphism without fixed points.
 %In particular, the two generators of $\Gamma$ fix exactly the two fibres coming from $B_i, B'_i$ in each elliptic fibration.
Hence it follows that the image of the elliptic fibration $\frakh_1$ by the map $\tilde \sigma^{-1}\tilde r$ in diagram (\ref{diagram2}) is an elliptic fibration on $ X'$. We will denote it by $\bar\frakh_1$.

Now statements $a)$, $b)$, $c)$ are easy consequences of the analogous statements in Proposition \ref{fib}.

According to Proposition \ref{hc} the cubics $B_1$ and $B'_1$ contain the $9$ vertices of the triangles $T_2, T_3, T_4$ in the Hesse pencil. 
Hence the fibres $E_1, E'_1$ contain each $9$ points in the preimage of the $9$ vertices by $q$.
It follows from Remark \ref{x/gamma} that the images of $E_1$ and $E'_1$
in $X/\Gamma$ contain each $3$ singular points of type $A_2$. 
The preimage of one of these fibres in the minimal resolution $X'$ is a fibre of type $IV^*$ in the elliptic fibration $\bar\frakh_1$ on $X'$ (the union of $3$ exceptional divisors of type $A_2$ and the proper transform of $E_1$ or $E'_1$). 
%We will denote these fibres by $\bar E_1$, $\bar E'_1$.

It can be easily seen that the $6$ singular fibres of $\frakh_1$ of type $IV$ belong to two orbits for the action of $\Gamma$.
In fact, the singular points in each of these fibres are the preimages by $q$ of the vertices of $T_1$ (see Proposition \ref{fib}).
The image of a singular fibre of type $IV$ in $X/\Gamma$ is a rational curve containing a singular point of type $A_2$ and its preimage in $X'$ is again a fibre of type $IV$. Hence $\bar\frakh_1$ has two fibres of type $IV$.
\qed
\begin{Rem}
%Recall that the fibres $F_1, F'_1$ are mapped by $q'$ to two cusps $q_1$, $q'_1$ of $C'_6$  (see section 6). 
%Let $\sigma$ be the covering involution of $q$. Since $B_1$ and $B'_1$ intersect exactly in the $9$ vertices of the triangles, then 
%$F_1$ and $ \sigma(F'_1)$ intersect in $9$ preimages of the vertices. Hence the three exceptional divisors of type $A_2$ in a fibre of type $IV^*$ intersect the image of the fibre $\sigma(F'_1)$ of $\frakh'_1$.
It can be proved that the image of any of these fibrations by the cover $q'$ is  a one-dimensional family of curves of degree $9$ in $\Ps 2$ with $8$ triple points in $q_1, \dots,q_8$ and $3$ cusps on $C_6'$.
In fact, $q'$ sends the fibre $\tilde{\sigma}^{-1}\tilde r (F_1)$ of $\bar\frakh_1$ to the union of the $6$ inflection lines through $q_1$ and $q'_1$ not containing $q_0$, where the $3$ lines through $q_1$ are double. Clearly, the similar statement is true for $F'_1$ (now the lines through $p'_1$ are double).
Hence the image of a fibre of $\bar\frakh_1$ is a plane curve $D$ of degree $9$ with $8$ triple points at $q_1,\dots,q_8$.
Moreover, the curve $D$ intersects the sextic $C'_6$ in $6$ more points and since its inverse image in $X'$ has genus one, then $D$ must also have three cusps at smooth points of $C'_6$ which are resolved in the double cover $q'$.
\end{Rem}

\begin{Thm}\label{kummer}
The $K3$ surface $X'$ is birationally isomorphic to the Kummer surface 
$$\Kum(E_{\e}\times E_{\e}) = (E_{\e}\times E_{\e})/(x \mapsto -x),$$ where $E_{\e}$ is the elliptic curve with fundamental periods $1,\e$. 
Its transcendental lattice has rank two and its intersection matrix with respect to a suitable basis is
$$A_2(-2)=\left(
\begin{array}{cc}
4&2\\
2& 4\end{array} \right 
). $$
 \end{Thm}

\begin{proof}
We will consider one of the jacobian fibrations on $X'$ described in Proposition \ref{iac}.
Let $M$ be the lattice generated by the two sections $e_1$, $e_2$, the components of the $3$ singular fibres of type $I_6$ not intersecting $e_2$ and the components of the fibre of type $I_3$. The intersection matrix of $M$ has determinant $-2^2\cdot 3^5$, hence $\rank\, M=\rank\, S_{X'}=20$ and $\rank\, T_{X'}=2$.

%Since $\Gamma$ acts symplectically on $X$, the transcendental lattice $T_{X'}$ has the same rank of $T_X$. 
The non-symplectic automorphism $g_4$ of order $3$ on $X$ induces an  automorphism $g_4'$ on $X'$. 
Recall that $g_4$ fixes the curve $R=\{x=0\}$ on $X$, hence $g_4'$  fixes the proper transform of $\tilde r(R)$ on $X'$.
Thus by Theorem \ref{auto} $g_4'$ is a non-symplectic automorphism of order three on $X'$.
This implies, as in the proof of Theorem \ref{k1}, that the intersection matrix of $T_{X'}$ is of the form \eqref{matrix} with respect to a proper choice of generators, in particular its discriminant group is isomorphic to $\Z/3\Z \oplus \Z/3m\Z $. 
%The projection rational map $X\to X'$ induces an embedding of the lattices $T_{X'}[9]\hookrightarrow T_{X}$, where $M[a]$ denoted the quadratic lattice whose quadratic form is multiplied by $a$. The known relationship between the discriminants of a lattice and its sublattice of finite index  $i$ gives 
%$$i^2\discr(T_{X}) = 2^2\cdot 3^3\cdot i^2 = \discr(T_X[9]) = 3^5\cdot m^2.$$ It is also known that the index $i$ divides the order of the group $\Gamma$, i.e. is equal to 1, 3, or 9. This gives the only possible value $m = 2$ and $i = 3$. 

A direct computation of $M^*$ shows that the discriminant group $A_M$ is isomorphic to $ \Z/3\Z^3\oplus \Z/6\Z^2$. Since $M$ is a sublattice of finite index of $S_{X'}$, the discriminant group $A_{T_{X'}}\cong A_{S_{X'}}$  is isomorphic to a quotient of a subgroup of $A_M$. This implies that $m\le 2$. 

By Theorem 3.1 (Table 1.1), \cite{N} the unique $K3$ surface with transcendental lattice as in \eqref{matrix} with $m=1$ has no jacobian elliptic fibration as in Proposition \ref{iac}.   Hence $m=2$ and by \cite{SI} $X'$ is isomorphic to the Kummer surface of the abelian surface $E_{\e}\times E_{\e}$.
\end{proof}
\begin{Rem}\ \\
i) 
In \cite{K} it is proved that all elliptic fibrations on the Kummer surface $\Kum(E_{\e}\times E_{\e})$ are jacobian. All these fibrations and their Mordell-Weil groups are described in \cite{N}. In particular it is proved that the Mordell-Weil group of the elliptic fibration in Proposition \ref{iac} is isomorphic to $\Z\oplus \Z/3\Z$ and those of
 the $8$ elliptic fibrations in Proposition \ref{fib2} are isomorphic to $\Z^2\oplus \Z/3\Z$ (see Theorem 3.1, Table 1.3, No. 19, 30). \\
ii) 
%By construction, we know that the quotient group $G_{216}/\Gamma\cong 2.A_4$ is a group of automorphisms of $X'$ and the action of its normal subgroup isomorphic to $Q_8$ is symplectic. 
The full automorphism group of $X'$ has been computed in \cite{KK} while the full automorphism group of $X$ is not known at present.

\end{Rem}
%\begin{Rem}
%Each fibre of type $A_2$ in the elliptic fibration $\mathcal F_i$ corresponds to a curve of degree nine given by: the union of the two lines $l,m$ containing $3$ cusps and not passing through $q_i$ and $q'_i$, a curve of degree $7$ passing through $l\cap m$ with $6$ nodes in the cusps on $l\cup m$,  triple points in $q_i$ and $q'_i$ and $3$ cusps on the sextic $C'$.
%\end{Rem}

\end{document}